\documentclass[a4,11pt]{article}
\usepackage{amsfonts,amssymb,amsthm,amsmath}
\usepackage[matrix,arrow]{xy}

\theoremstyle{plain}

\newtheorem*{Th*}{Theorem}

\newtheorem*{Cor*}{Corollary}

\theoremstyle{definition}

\theoremstyle{remark}

\numberwithin{equation}{section}

\makeatletter
\def\Set@Scallop[#1]#2#3{{#1}\Parens{#2}{#3}}
\newcommand\DeclareScalableOperator[2]{%
  \expandafter\def\csname#1\endcsname{\@ifnextchar[{{#2}\Set@Scallop}{{#2}\Set@Scallop[{}]}}
}
\makeatother

\DeclareScalableOperator{Ct}{\mathcal{C}} % continuous functions/maps
\DeclareScalableOperator{Cc}{\mathcal{K}} % continuous w/ cpt support
\DeclareScalableOperator{Lp}{\mathbf{L}} % L^p spaces
\DeclareScalableOperator{Hom}{\mathrm{Hom}} % homomorphisms
\DeclareScalableOperator{End}{\mathrm{End}} 
\DeclareScalableOperator{Aut}{\mathrm{Aut}} 
\DeclareScalableOperator{Uenv}{\mathfrak{U}} % universal enveloping algebra

\newcommand{\Fa}{For all }
\newcommand{\fa}{for all }

\newcommand{\fs}{for some }
\newcommand\mathfa[1][{}]{\quad\text{\fa{#1} }}

\newcommand{\scth}{such that }
\newcommand{\AND}{and}

\newcommand\mathtxt[1]{\quad\text{{#1}}\quad}

\newcommand{\nd}{\mathtxt\AND}
\newcommand{\nda}{\ \text{\AND}\ }

\newcommand\eps{\varepsilon}
\newcommand\nats{\mathbb{N}}
\newcommand\ints{\mathbb{Z}}

\newcommand\reals{\mathbb{R}}
\newcommand\cplxs{\mathbb{C}}

\newcommand\vvoid{\varnothing}
\newcommand\sle{\leqslant}
\newcommand\sge{\geqslant}

\DeclareMathOperator\dom{\mathrm{dom}}

\DeclareMathOperator\Ad{\mathrm{Ad}}
\DeclareMathOperator\ad{\mathrm{ad}}
\DeclareMathOperator\GL{\mathrm{GL}}

\DeclareMathOperator\sll{\ger{sl}}

\makeatletter
\newcommand\Size[7][1]{% produces #3#4#5#6#7, where #3,#5,#7 are operators of size #2(0-4); if #1 is 0, #5 is not resized
                                 \ifx#20%
                                        \def\r@l{}\def\r@m{}\def\r@r{}%
                                 \else%
                                    \ifx#21%
                                           \def\r@l{\bigl}\def\r@r{\bigr}\def\r@m{\bigm}%
                                    \else%
                                           \ifx#22%
                                                 \def\r@l{\Bigl}\def\r@r{\Bigr}\def\r@m{\Bigm}%
                                            \else%
                                                 \ifx#23%
                                                        \def\r@l{\biggl}\def\r@r{\biggr}\def\r@m{\biggm}%
                                                  \else
                                                        \ifx#24%
                                                        \def\r@l{\Biggl}\def\r@r{\Biggr}\def\r@m{\Biggm}%
                                                        \fi%
                                                  \fi%
                                            \fi%
                                      \fi%
                                 \fi%
                                 \ifx#10%
                                       \def\r@m{}%
                                 \fi%
                                 \r@l#3{#4}\r@m#5{#6}\r@r#7%
}%

\makeatother

\newcommand\Set[3]{% set theoretical brackets #1:size (0-4), #2:elements, #3:defining predicate
                                 \Size{#1}{\{}{#2}{|}{#3}{\}}%
}%
\newcommand\Cdual[3]{% Bilinear inner product of size #1 (0-4)
                                 \Size[0]{#1}{\langle}{#2}{,}{#3}{\rangle}%
}%
\newcommand\Parens[2]{% Parentheses of size #1 (0-4)
  \Size[0]{#1}{(}{#2}{}{}{)}
}

\newcommand\Abs[2]{% absolute value of size #1 (0-4)
  \Size[0]{#1}{\lvert}{#2}{}{}{\rvert}
}

\makeatletter

\newcommand{\IfUpperCase}[1]{\begingroup 
  \protected@edef\@tempa{\expandafter\@firstofone\@firstofone#1.}%
  \expandafter\IfUpperCasE \@tempa\delimiter}

\def\IfUpperCasE #1#2\delimiter{%
  \protected@edef\@tempa{\meaning#1\meaning a}%
  \ifnum \expandafter\IfUppercaSE\@tempa \IfUppercaSE
   \endgroup \expandafter\@firstoftwo
  \else
   \endgroup \expandafter\@secondoftwo
  \fi}

\@ifundefined{strip@prefix}{\def\strip@prefix#1>{}}{}
\def\@tempa{the letter }
\edef\@tempa{\expandafter\strip@prefix\meaning\@tempa}

\expandafter\def\expandafter\IfUppercaSE\expandafter#\expandafter1\@tempa#2#3\IfUppercaSE{\uccode`#2=`#2 }

\newif\ifuc@se
\def\setuc@se#1{\IfUpperCase{#1}{\uc@setrue}{\uc@sefalse}}

\def\theoremn@me#1{\ifuc@se \lowercase{\csname#1name\endcsname}\ignorespaces%
  \else \edef\@temp{\lowercase{\lowercase{\csname#1name\endcsname}}}\@temp\ignorespaces%
  \fi}
\def\theoremn@mes#1{\ifuc@se \lowercase{\csname#1names\endcsname}\ignorespaces%
  \else \edef\@temp{\lowercase{\lowercase{\csname#1names\endcsname}}}\@temp\ignorespaces%
  \fi}
  
\def\thmref#1#2{\setuc@se{#1}\lowercase{{\theoremn@me{#1}\lowercase{\ref{#1:#2}}}}}

\newcommand{\DefTheorem}[2]{\newenvironment{#1}[1][\empty]{\ignorespaces\begin{#2}\ifx##1\empty{}\else\lowercase{\label{#1:##1}}\fi\ignorespaces}{\end{#2}\ignorespacesafterend}}

\DefTheorem{Th}{theorem}
\DefTheorem{Prop}{proposition}
\DefTheorem{Cor}{corollary}
\DefTheorem{Lem}{lemma}
\DefTheorem{Def}{definition}
\DefTheorem{Rem}{remark}
\DefTheorem{Par}{para}

\newenvironment{Par*}{\ignorespaces\noindent\ignorespaces}{\ignorespacesafterend}
\makeatletter

\newif\if@smallmat
\newif\if@none
\newif\if@paren
\newif\if@brack
\newif\if@brace
\newif\if@vline
\newenvironment{Matrix}[2][1]% matrices #1:size (0-1) #2:paren type (0-5, none, parentheses, brackets, braces, vertical lines, vertical double lines)
                                 {\ifx#20%
                                        \@smallmattrue%
                                  \else%
                                         \@smallmatfalse
                                  \fi%
                                  \ifx#11%
                                         \@nonefalse\@parentrue\@brackfalse\@bracefalse\@vlinefalse%
                                  \else%
                                       \ifx#12%
                                            \@nonefalse\@parenfalse\@bracktrue\@bracefalse\@vlinefalse%
                                        \else%
                                            \ifx#13%
                                                 \@nonefalse\@parenfalse\@brackfalse\@bracetrue\@vlinefalse%
                                            \else%
                                                 \ifx#14%
                                                       \@nonefalse\@parenfalse\@brackfalse\@bracefalse\@vlinetrue
                                                 \else%
                                                       \ifx#15%
                                                             \@nonefalse\@parenfalse\@brackfalse\@bracefalse\@vlinefalse%
                                                       \else%
                                                             \@nonetrue\@parenfalse\@brackfalse\@bracefalse\@vlinefalse%
                                                       \fi%
                                                 \fi%
                                            \fi%
                                        \fi%
                                   \fi%
                                   \if@smallmat%
                                        \if@none%
                                             \begin{smallmatrix}%
                                        \else%
                                            \if@paren%
                                                  \bigl(\begin{smallmatrix}%
                                            \else%
                                                  \if@brack%
                                                          \bigl[\begin{smallmatrix}%
                                                  \else%
                                                          \if@brace%
                                                               \bigl\{\begin{smallmatrix}%
                                                          \else%
                                                               \if@vline%
                                                                    \bigl\lvert\begin{smallmatrix}%
                                                                \else%
                                                                    \bigl\lVert\begin{smallmatrix}%
                                                                \fi%
                                                          \fi%
                                                  \fi%
                                            \fi%
                                        \fi%
                                   \else%
                                        \if@none%
                                             \begin{matrix}%
                                        \else%
                                            \if@paren%
                                                  \begin{pmatrix}%
                                            \else%
                                                  \if@brack%
                                                          \begin{bmatrix}%
                                                  \else%
                                                          \if@brace%
                                                               \begin{Bmatrix}%
                                                          \else%
                                                               \if@vline%
                                                                    \begin{vmatrix}%
                                                                \else%
                                                                    \begin{Vmatrix}%
                                                                \fi%
                                                          \fi%
                                                  \fi%
                                            \fi%
                                        \fi%
                                   \fi}%
                                  {\if@smallmat%
                                        \if@none%
                                             \end{smallmatrix}%
                                        \else%
                                            \if@paren%
                                                  \end{smallmatrix}\bigr)%
                                            \else%
                                                  \if@brack%
                                                          \end{smallmatrix}\bigr]%
                                                  \else%
                                                          \if@brace%
                                                               \end{smallmatrix}\bigr\}%
                                                          \else%
                                                               \if@vline%
                                                                    \end{smallmatrix}\bigr\rvert%
                                                                \else%
                                                                    \end{smallmatrix}\bigr\rVert%
                                                                \fi%
                                                          \fi%
                                                  \fi%
                                            \fi%
                                         \fi%
                                   \else%
                                        \if@none%
                                             \end{matrix}%
                                        \else%
                                            \if@paren%
                                                  \end{pmatrix}%
                                            \else%
                                                  \if@brack%
                                                          \end{bmatrix}%
                                                  \else%
                                                          \if@brace%
                                                               \end{Bmatrix}%
                                                          \else%
                                                               \if@vline%
                                                                    \end{vmatrix}%
                                                                \else%
                                                                    \end{Vmatrix}%
                                                                \fi%
                                                          \fi%
                                                  \fi%
                                            \fi%
                                        \fi%
                                   \fi}%

\makeatother

\def\ger{\mathfrak}
\DeclareMathOperator\str{\mathrm{str}}

\title{Chevalley's restriction theorem\\ for reductive symmetric superpairs}
\author{Alexander Alldridge, Joachim Hilgert, and Martin R.~Zirnbauer}

\CompileMatrices

\begin{document}

\maketitle

\begin{abstract}
	Let $(\ger g,\ger k)$ be a reductive symmetric superpair of even type, \emph{i.e.}~so that there exists an even Cartan subspace $\ger a\subset\ger p$. The restriction map $S(\ger p^*)^{\ger k}\to S(\ger a^*)^W$ where $W=W(\ger g_0:\ger a)$ is the Weyl group, is injective. We determine its image explicitly. 
	
	In particular, our theorem applies to the case of a symmetric superpair of group type, \emph{i.e.}~$(\ger k\oplus\ger k,\ger k)$ with the flip involution where $\ger k$ is a classical Lie superalgebra with a non-degenerate invariant even form (equivalently, a finite-dimensional contragredient Lie superalgebra). Thus, we obtain a new proof of the generalisation of Chevalley's restriction theorem due to Sergeev and Kac, Gorelik.
	
	For general symmetric superpairs, the invariants exhibit a new and surprising behaviour. We illustrate this phenomenon by a detailed discussion in the example $\ger g=C(q+1)=\ger{osp}(2|2q,\cplxs)$, endowed with a special involution. Here, the invariant algebra defines a singular algebraic curve. 
\end{abstract}

\section{Introduction}

The physical motivation for the development of supermanifolds stems from quantum field theory in its functional integral formulation, which describes fermionic particles by anticommuting fields. In the 1970s, pioneering work by Berezin strongly suggested that commuting and anticommuting variables should be treated on equal footing. Several theories of supermanifolds have been advocated, among which the definition of Berezin, Kostant, and Leites is one of the most commonly used in mathematics. 

Our motivation for the study of supermanifolds comes from the study of certain nonlinear $\sigma$-models with supersymmetry. Indeed, it is known from the work of the third named author \cite{zirnbauer-rsss} that Riemannian symmetric superspaces occur naturally in the large $N$ limit of certain random matrix ensembles, which correspond to Cartan's ten infinite series of symmetric spaces. In spite of their importance in physics, the mathematical theory of these superspaces is virtually non-existent. (But compare \cite{duflo-petracci-formalss,zirnbauer-superbos,goertsches-diss}.) We intend to initiate the systematic study of Riemannian symmetric superspaces, in order to obtain a good understanding of, in particular, the invariant differential operators, the spherical functions, and the related harmonic analysis. 
The present work lays an important foundation for this endeavour: the generalisation of Chevalley's restriction theorem to the super setting. 

To describe our results in detail, let us make our assumptions more precise. Let $\ger g$ be a complex Lie superalgebra with even centre \scth $\ger g_0$ is reductive in $\ger g$ and $\ger g$ carries an even invariant supersymmetric form. Let $\theta$ be an involutive automorphism of $\ger g$, and denote by $\ger g=\ger k\oplus\ger p$ the decomposition into $\theta$-eigenspaces. We say that $(\ger g,\ger k)$ is a \emph{reductive superpair}, and it is of \emph{even type} if there exists an even Cartan subspace $\ger a\subset\ger p_0$. 

Assume that $(\ger g,\ger k)$ is a reductive symmetric superpair of even type. Let $\bar\Sigma_1^+$ denote the set of positive roots of $\ger g_1:\ger a$ \scth $\lambda,2\lambda$ are no roots of $\ger g_0:\ger a$. To each $\lambda\in\bar\Sigma_1^+$, one associates a set $\mathcal R_\lambda$ of differential operators with rational coefficients on $\ger a$. 

Our main results are as follows. 

\begin{Th*}[A]
	Let $I(\ger a^*)$ be the image of the restriction map $S(\ger p^*)^{\ger k}\to S(\ger a^*)$ (which is injective). Then $I(\ger a^*)$ is the set of $W$-invariant polynomials on $\ger a$ which lie in the common domain of all operators in $\mathcal R_\lambda$, $\lambda\in\bar\Sigma_1^+$. Here, $W$ is the Weyl group of $\ger g_0:\ger a$.
\end{Th*}

For $\lambda\in\bar\Sigma_1^+$, let $A_\lambda\in\ger a$ be the corresponding coroot, and denote by $\partial(A_\lambda)$ the directional derivative operator in the direction of $A_\lambda$. Then the image $I(\ger a^*)$ can be characterised in more explicit terms, as follows.

\begin{Th*}[B]
	We have $I(\ger a^*)=\bigcap_{\lambda\in\bar\Sigma_1^+}S(\ger a^*)^W\cap I_\lambda$ where
	\[
		I_\lambda=\textstyle\bigcap_{j=1}^{\frac12m_{1,\lambda}}\dom\lambda^{-j}\partial(A_\lambda)^j\mathtxt{if}\lambda(A_\lambda)=0\ ,
	\]
	and if $\lambda(A_\lambda)\neq0$, then $I_\lambda$ consists of those $p\in\cplxs[\ger a]$ \scth
	\[
		\partial(A_\lambda)^kp|_{\ker\lambda}=0\mathfa[odd integers] k\ ,\ 1\sle k\sle m_{1,\lambda}-1\ .
	\]
	Here, $m_{1,\lambda}$ denotes the multiplicity of $\lambda$ in $\ger g_1$ (and is an even integer).
\end{Th*}

If the symmetric pair $(\ger g,\ger k)$ is of \emph{group type}, \emph{i.e.}~$\ger g=\ger k\oplus\ger k$ with the flip involution, then \fa $\lambda\in\bar\Sigma_1^+$, $\lambda(A_\lambda)=0$, and the multiplicity $m_{1,\lambda}=2$. In this case, Theorem (B) reduces to $I(\ger a^*)=\bigcap_{\lambda\in\bar\Sigma_1^+}S(\ger a^*)^W\cap\dom\lambda^{-1}\partial(A_\lambda)$. The situation where $\lambda(A_\lambda)\neq0$ \fs $\lambda\in\bar\Sigma_1^+$ occurs if and only if $\ger g$ contains symmetric subalgebras $\ger s\cong C(2)=\ger{osp}(2|2)$ where $\ger s_0\cap\ger k=\sll(2,\cplxs)$. This is case for $\ger g=C(q+1)$ with a special involution, and in this case, the invariant algebra $I(\ger a^*)$ defines the singular curve $z^{2q+1}=w^2$ (\thmref{Cor}{singular}). 

Let us place our result in the context of the literature. The Theorems (A) and (B) apply to the case of classical Lie superalgebras with non-degenerate invariant even form (equivalently, finite-dimensional contragredient Lie superalgebras), considered as symmetric superspaces of group type. In this case, the result is due to Sergeev \cite{sergeev-invpol}, Kac \cite{kac-laplace}, and Gorelik \cite{gorelik-kacconstruction}, and we simply furnish a new (and elementary) proof. (The results of Sergeev are also valid for basic Lie superalgebras which are not contragredient.) For some particular cases, there are earlier results by Berezin \cite{berezin}.

Sergeev's original proof involves case-by-case calculations. The proof by Gorelik---which carries out in detail ideas due to Kac in the context of Kac--Moody algebras---is classification-free, and uses so-called Shapovalov determinants. Moreover, the result of Kac and Gorelik actually characterises the image of the Harish-Chandra homomorphism rather than the image of the restriction map on the symmetric algebra, and is therefore more fundamental than our result.

Still in the case of symmetric superpairs of group type, Kac \cite{kac-typical} and Santos \cite{santos-zuckerman} describe the image of the restriction morphism in terms of supercharacters of certain (cohomologically) induced modules (instead of a characterisation in terms of a system of differential equations). This approach cannot carry over to the case of symmetric pairs, as is known in the even case from the work of Helgason \cite{helgason-fundamental}. 

Our result also applies in the context of Riemannian symmetric superspaces, where one has an even non-degenerate $\mathcal G$-invariant supersymmetric form on $\mathcal G/\mathcal K$ whose restriction to the base $G/K$ is Riemannian. In this setting, it is to our knowledge completely new and not covered by earlier results. We point out that a particular case was proved in the PhD thesis of Fuchs \cite{fuchs-diss}, in the framework of the `supermatrix model', using a technique due to Berezin. 

In the context of harmonic analysis of even Riemannian symmetric spaces $G/K$, Chevalley's restriction theorem enters crucially, since it determines the image of the Harish-Chandra homomorphism, and thereby, the spectrum of the algebra $\mathbb D(G/K)$ of $G$-invariant differential operators on $G/K$. It is an important ingredient in the proof of Harish-Chandra's integral formula for the spherical functions. In a series of forthcoming papers, we will apply our generalisation of Chevalley's restriction theorem to obtain analogous results in the context of Riemannian symmetric superspaces. 

\medskip\noindent
Let us give a brief overview of the contents of our paper. We review some basic facts on root decompositions in sections 2.1-2.2. In section 2.3, we introduce our main tool in the proof of Theorem (A), a certain twisted action $u_z$ on the supersymmetric algebra $S(\ger p)$. In section 3.1, we define the `radial component' map $\gamma_z$ via the twisted action $u_z$. The proofs of Theorems (A) and (B) are contained in sections 3.2 and 3.3, respectively. The former comes down to a study of the singularities of $\gamma_z$ as a function of the semi-simple $z\in\ger p_0$, whereas the latter consists in an elementary and explicit discussion of the radial components of certain differential operators. In sections 4.1 and 4.2, we discuss the generality of the `even type' condition, and study an extreme example in some detail. 

The first named author wishes to thank C.~Torossian (Paris VII) for his enlightening remarks on a talk given on an earlier version of this paper. The first and second named author wish to thank M.~Duflo (Paris VII) for helpful discussions, comments, and references. The second named author wishes to thank K.~Nishiyama (Kyoto) for several discussions on the topic. Last, not least, we wish to thank an anonymous but diligent referee whose suggestions greatly improved the presentation of our main technical device. 

This research was partly funded by the IRTG ``Geometry and Analysis of Symmetries'', supported by Deutsche Forschungsgemeinschaft (DFG), Minist\`ere de l'\'Education Nationale (MENESR), and Deutsch-Franz\"osische Hochschule (DFH-UFA), and by the SFB/Transregio 12 ``Symmetry and Universality in Mesoscopic Systems'', supported by Deutsche Forschungsgemeinschaft (DFG).

\section{Some basic facts and definitions}

\begin{Par*}
	In this section, we mostly collect some basic facts concerning (restricted) root decompositions of Lie superalgebras, and the (super-) symmetric algebra, along with some definitions which we find useful to formulate our main results. As general references for matters super, we refer the reader to \cite{kostant-supergeom,deligne-morgan-susy,kac-liesuperalgs,scheunert-liesuperalgs}
\end{Par*}

\subsection{Roots of a basic quadratic Lie superalgebra}

\begin{Def}
	Let $\ger g=\ger g_0\oplus\ger g_1$ be a Lie superalgebra over $\cplxs$ and $b$ a bilinear form $b$. Recall that $b$ is \emph{supersymmetric} if $b(u,v)=(-1)^{\Abs0u\Abs0v}b(v,u)$ \fa homogeneous $u,v$. We shall call $(\ger g,b)$ \emph{quadratic} if $b$ is a non-degenerate, $\ger g$-invariant, even and supersymmetric form on $\ger g$. We shall say that $\ger g$ is \emph{basic} if $\ger g_0$ is reductive in $\ger g$ (\emph{i.e.}~$\ger g$ is a semi-simple $\ger g_0$-module) and $\ger z(\ger g)\subset\ger g_0$ where $\ger z(\ger g)$ denotes the centre of $\ger g$.
\end{Def}

\begin{Par}
	Let $(\ger g,b)$ be a basic quadratic Lie superalgebra, and $\ger b$ be a Cartan subalgebra of $\ger g_0$. 
	
	As usual \cite[Chapter II, \S~4.6]{scheunert-liesuperalgs}, we define
  \[
  V^\alpha=\Set1{x\in V}{\exists\,n\in\nats\,:\,(h-\alpha(h))^n(x)=0\ \text\fa h\in\ger b}\quad,\quad\alpha\in\ger b^*
  \]
  for any $\ger b$-module $V$. Further, the sets of even resp.~odd roots for $\ger b$ are  
  \[
  \Delta_0(\ger g:\ger b)
  =\Set1{\alpha\in\ger b^*\setminus0}{\ger g_0^\alpha\neq0}
  \nd\Delta_1(\ger g:\ger b)
  =\Set1{\alpha\in\ger b^*}{\ger g_1^\alpha\neq0}\ .
  \]
  We also write $\Delta_j=\Delta_j(\ger g:\ger b)$. Let $\Delta=\Delta(\ger g:\ger b)=\Delta_0\cup\Delta_1$. The elements of $\Delta$ are called \emph{roots}. We have 
  \[
  \ger g=\ger b\oplus\textstyle\bigoplus_{\alpha\in\Delta}\ger g^\alpha
  =\ger b\oplus\bigoplus_{\alpha\in\Delta_0}\ger g_0^\alpha
  	\oplus\bigoplus_{\alpha\in\Delta_1}\ger g_1^\alpha\ .
  \]
  It is obvious that $\Delta_0=\Delta(\ger g_0:\ger b)$, so in particular, it is a reduced abstract root system in its real linear span. Also, since $\ger g_0$ is reductive in $\ger g$, the root spaces $\ger g_i^\alpha$ are the joint eigenspaces of $\ad h$, $h\in\ger b$ (and not only generalised ones).
  
  We collect the basic statements about $\ger b$-roots. The results are known (\emph{e.g.}~\cite{scheunert-liesuperalgs,benayadi-root}), so we omit their proofs.
\end{Par}

\begin{Prop}[superroot]
Let $\ger g$ be a basic quadratic Lie superalgebra with invariant form $b$, and $\ger b$ a Cartan subalgebra of $\ger g_0$.
\begin{enumerate}
\item For $\alpha,\beta\in\Delta\cup0$, we have $b(\ger g_j^\alpha,\ger g_k^\beta)=0$ unless $j=k$ and $\alpha=-\beta$.
\item The form $b$ induces a non-degenerate pairing $\ger g_j^\alpha\times\ger g_j^{-\alpha}\to\cplxs$. In particular, we have $\dim\ger g_j^\alpha=\dim\ger g_j^{-\alpha}$ and $\Delta_j=-\Delta_j$ for $j\in\ints/2\ints$. 
\item The form $b$ is non-degenerate on $\ger b$, so for any $\lambda\in\ger b^*$, there exists a unique $h_\lambda\in\ger b$ \scth $b(h_\lambda,h)=\lambda(h)$ \fa $h\in\ger b$. 
\item If $\alpha(h_\alpha)\neq0$, $\alpha\in\Delta_1$, then $2\alpha\in\Delta_0$. In particular, $\Delta_0\cap\Delta_1=\vvoid$. 
\item We have $\ger g_1^0=\ger z_1(\ger g)=\{x\in\ger g_1\mid[x,\ger g]=0\}=0$, so $0\not\in\Delta_1$.
\item All root spaces $\ger g^\alpha$, $\alpha\in\Delta$, $\alpha(h_\alpha)\neq0$, are one-dimensional. 
\end{enumerate}
\end{Prop}

\subsection{Restricted roots of a reductive symmetric superpair}

\begin{Def}
Let $(\ger g,b)$ be a complex quadratic Lie superalgebra, and $\theta:\ger g\to\ger g$ an involutive automorphism leaving the form $b$ invariant. If $\ger g=\ger k\oplus\ger p$ is the $\theta$-eigenspace decomposition, then we shall call $(\ger g,\ger k)$ a \emph{symmetric superpair}.  We shall say that $(\ger g,\ger k)$ is \emph{reductive} if, moreover, $\ger g$ is basic.

Note that for any symmetric superpair $(\ger g,\ger k)$, $\ger k$ and $\ger p$ are $b$-orthogonal and non-degenerate. It is also useful to consider the form $b^\theta(x,y)=b(x,\theta y)$ which is even, supersymmetric, non-degenerate and $\ger k$-invariant.

Let $(\ger g,\ger k)$ be a reductive symmetric superpair. For arbitrary subspaces $\ger c,\ger d\subset\ger g$, let $\ger z_{\ger d}(\ger c)=\Set0{d\in\ger d}{[d,\ger c]=0}$ denote the centraliser of $\ger c$ in $\ger d$. Any linear subspace $\ger a=\ger z_{\ger p}(\ger a)\subset\ger p_0$ consisting of semi-simple elements of $\ger g_0$ is called an \emph{even Cartan subspace}. If an even Cartan subspace exists, then we say that $(\ger g,\ger k)$ is of \emph{even type}. 
\end{Def}

\begin{Par*}
	We state some generalities on even Cartan subspaces. These are known and straightforward to deduce from standard texts such as \cite{dixmier-envalg,borel-rss}.
\end{Par*}

\begin{Lem}
Let $\ger a\subset\ger g$ be an even Cartan subspace.
\begin{enumerate}
\item $\ger a$ is reductive in $\ger g$, \emph{i.e.}~$\ger g$ is a semi-simple $\ger a$-module.
\item $\ger z_{\ger g_0}(\ger a)$ and $\ger z_{\ger g_1}(\ger a)$ are $b$-non-degenerate.
\item $\ger z_{\ger g_0}(\ger a)=\ger m_0\oplus\ger a$ and $\ger z_{\ger g_1}(\ger a)=\ger m_1$ where $\ger m_i=\ger z_{\ger k_i}(\ger a)$, and the sum is $b$-orthogonal.
\item $\ger m_0$, $\ger m_1$, and $\ger a$ are $b$-non-degenerate.  
\item There exists a $\theta$-stable Cartan subalgebra $\ger b$ of $\ger g_0$ containing $\ger a$. 
\end{enumerate}
\end{Lem}

\begin{Par}[grouptype]
	Let $\ger k$ be a classical Lie superalgebra with a non-degenerate invariant even form $B$ \cite{kac-repnclassical}. Then $\ger k_0$ is reductive in $\ger k$, and $\ger z(\ger k)$ is even. We may define $\ger g=\ger k\oplus\ger k$, and $b(x,y,x',y')=B(x,x')+B(y,y')$. Then $(\ger g,b)$ is basic quadratic. The flip involution $\theta(x,y)=(y,x)$ turns $(\ger g,\ger k)$ into a reductive symmetric superpair (where $\ger k$ is, as is customary, identified with the diagonal in $\ger g$). We call such a pair of \emph{group type}. 
	
	Moreover, any Cartan subalgebra $\ger a$ of $\ger k_0$ yields an even Cartan subspace for the superpair $(\ger g,\ger k)$. Indeed, $\ger p=\Set1{(x,-x)}{x\in\ger k}$, and the assertion follows from \thmref{Prop}{superroot}~(v).
\end{Par}

\begin{Par}
	In what follows, let $(\ger g,\ger k)$ be a reductive symmetric superpair of even type, $\ger a\subset\ger p$ an even Cartan subspace, and $\ger b\subset\ger g_0$ a $\theta$-stable Cartan subalgebra containing $\ger a$. The involution $\theta$ acts on $\ger b^*$ by $\theta\alpha=\alpha\circ\theta$ \fa $\alpha\in\ger b^*$. Let $\alpha_\pm=\frac12(1\pm\theta)\alpha$ \fa $\alpha\in\ger b^*$, and set
	\[
		\Sigma_j=\Sigma_j(\ger g:\ger a)=\Set1{\alpha_-}{\alpha\in\Delta_j\ ,\ \alpha\neq\theta\alpha}\ ,\ \Sigma=\Sigma(\ger g:\ger a)=\Sigma_0\cup\Sigma_1\ .
	\]
	(The union might not be disjoint.) Identifying $\ger a^*$ with the annihilator of $\ger b\cap\ger k$ in $\ger b^*$, these may be considered as subsets of $\ger a^*$. The elements of $\Sigma_0$, $\Sigma_1$, and $\Sigma$ are called \emph{even restricted roots}, \emph{odd restricted roots}, and \emph{restricted roots}, respectively. For $\lambda\in\Sigma$, let 
	\[
		\Sigma_j(\lambda)=\Set1{\alpha\in\Delta_j}{\lambda=\alpha_-}\ ,\ \Sigma(\lambda)=\Sigma_0(\lambda)\cup\Sigma_1(\lambda)\ .
	\]
	In the following lemma, observe that $\lambda\in\Sigma_j(\lambda)$ means that $\lambda\in\Delta_j$. We omit the simple proof, which is exactly the same as in the even case \cite[Chapter 1.1, Appendix 2, Lemma 1]{warner-vol1}.
\end{Par}

\begin{Lem}[rootrestr-fibre-even]
	Let $\lambda\in\Sigma_j$, $j=0,1$. The map $\alpha\mapsto-\theta\alpha$ is a fixed point free involution of $\Sigma_j(\lambda)\setminus\lambda$. In particular, the cardinality of this set is even. 
\end{Lem}

\begin{Par}
	For $\lambda\in\Sigma$, let 
	\[
		\ger g_{j,\ger a}^\lambda=\Set1{x\in\ger g_j}{\forall h\in\ger a\,:\,[h,x]=\lambda(h)\cdot x}\ ,\ \ger g_{\ger a}^\lambda=\ger g_{0,\ger a}^\lambda\oplus\ger g_{1,\ger a}^\lambda\ ,
	\]
	and $m_{j,\lambda}=\dim_\cplxs\ger g^\lambda_{j,\ger a}$, the \emph{even} or \emph{odd multiplicity} of $\lambda$, according to whether $j=0$ or $j=1$. It is clear that
	\[
		\textstyle\ger g_{j,\ger a}^\lambda=\bigoplus_{\alpha\in\Sigma_j(\lambda)}\ger g^\alpha_j\ ,\ m_{j,\lambda}=\sum_{\alpha\in\Sigma_j(\lambda)}\dim_\cplxs\ger g_j^\alpha\ ,\nda\ger g=\ger z_{\ger g}(\ger a)\oplus\bigoplus_{\lambda\in\Sigma}\ger g^\lambda_{\ger a}\ .
	\]
\end{Par}

\begin{Par*}
	The following facts are certainly well-known. Lacking a reference, we give the short proof. 
\end{Par*}

\begin{Prop}[res-superroot]
Let $\alpha,\beta\in\Delta$, $\lambda\in\Sigma$, and $j,k\in\{0,1\}$.
\begin{enumerate}
	\item The form $b^\theta$ is zero on $\ger g^\alpha_j\times\ger g_k^\beta$, unless $j=k$ and $\alpha=-\theta\beta$, in which case it gives a non-degenerate pairing.
	\item There exists a unique $A_\lambda\in\ger a$ \scth $b(A_\lambda,h)=\lambda(h)$ \fa $h\in\ger a$. 
	\item We have $\dim_\cplxs\ger g_j^\alpha=\dim_\cplxs\ger g_j^{-\theta\alpha}$. 
	\item The subspace $\ger g_j(\lambda)=\ger g_{j,\ger a}^\lambda\oplus\ger g_{j,\ger a}^{-\lambda}$ is $\theta$-invariant and decomposes into $\theta$-eigenspaces as $\ger g_j(\lambda)=\ger k_j^\lambda\oplus\ger p_j^\lambda$.
	\item The odd multiplicity $m_{1,\lambda}$ is even, and $b^\theta$ defines a symplectic form on both $\ger k_1^\lambda$ and $\ger p_1^\lambda$. 
\end{enumerate}
\end{Prop}

\begin{proof}
	The form $b^\theta$ is even, so $b^\theta(\ger g_0,\ger g_1)=0$. For $x\in\ger g^\alpha_j$, $y\in\ger g^\beta_j$, we compute, \fa $h\in\ger b$, 
	\begin{align*}
		(\alpha+\theta\beta)(h)b^\theta(x,y)&=b^\theta([h,x],y)+b^\theta(x,[\theta h,y])\\
		&=b^\theta([h,x]+[x,h],y)=0\ .
	\end{align*}
	Hence, $b^\theta(x,y)=0$ if $\alpha\neq-\theta\beta$. Since $b^\theta$ is non-degenerate and $\ger g/\ger b$ is the sum of root spaces, $b^\theta$ induces a non-degenerate pairing of $\ger g_j^\alpha$ and $\ger g_j^{-\theta\alpha}$. We also know already that $\ger a$ is non-degenerate for $b^\theta$, and (i)-(iii) follow. Statement (iv) is immediate. 
	
	We have 
	\[
		\ger g_{1,\ger a}^\lambda/\ger g_1^\lambda\cong\textstyle\bigoplus_{\alpha\in\Sigma_j(\lambda)\setminus\lambda}\ger g_1^\alpha\ .
	\]
	By (iii) and \thmref{Lem}{rootrestr-fibre-even}, this space is even-dimensional. But $\lambda$ is a root if and only if $\lambda=-\theta\lambda$. Then $b^\theta$ defines a symplectic form on $\ger g_1^\lambda$ by (i), and this space is even-dimensional. Thus, $m_{1,\lambda}$ is even, and again by (i), $\ger g_{1,\ger a}^\lambda$ is $b^\theta$-non-degenerate. It is clear that $\ger k_1^\lambda$ and $\ger p_1^{\lambda}$ are $b^\theta$-non-degenerate because $\ger g_{1,\ger a}^\lambda$ and $\ger g_{1,\ger a}^{-\lambda}$ are. Hence, we obtain assertion (v).
\end{proof}

\begin{Rem}
	Unlike the case of unrestricted roots, there may exist $\lambda\in\Sigma_1$ \scth $2\lambda\not\in\Sigma$ but $\lambda$ is still anisotropic, \emph{i.e.}~$\lambda(A_\lambda)\neq0$. Indeed, consider $\ger g=\ger{osp}(2|2,\cplxs)$ ($\cong\sll(2|1,\cplxs)$). Then $\ger g_0=\ger o(2,\cplxs)\oplus\ger{sp}(2,\cplxs)=\ger{gl}(2,\cplxs)$ and $\ger g_1$ is the sum of the fundamental representation of $\ger g_0$ and its dual. 
	
	Define the involution $\theta$ to be conjugation by the element $\begin{Matrix}0\sigma&0\\0&1_2\end{Matrix}$ where $\sigma=\begin{Matrix}00&1\\1&0\end{Matrix}$. One finds that $\ger k_0=\sll(2,\cplxs)$ and $\ger p_0=\ger a=\ger z(\ger g_0)$ which is one-dimensional and non-degenerate for the supertrace form $b$. On the other hand, $\ger g_1=\ger g_1(\lambda)$ is the sum of the root spaces for certain odd roots $\pm\alpha$, $\pm\theta\alpha$ which restrict to $\pm\lambda$. Clearly, there are no even roots, so $2\lambda$ is not a restricted root. Since $A_\lambda$ generates $\ger a$, it is a $b$-anisotropic vector. We discuss this issue at some length in section 4.2. 
	
	We point out that it is also not hard to prove that any such root $\lambda$ occurs in this setup. \emph{I.e.}, given a reductive symmetric superpair $(\ger g,\ger k)$, for any $\lambda\in\Sigma_1$, $2\lambda\not\in\Sigma$, $\lambda(A_\lambda)\neq0$, there exists a $b$-non-degenerate $\theta$-invariant subalgebra $\ger s\cong\ger{osp}(2|2,\cplxs)$ \scth $\ger p\cap\ger s_0=\cplxs A_\lambda=\ger z(\ger s_0)$ (the centre of $\ger s_0$), and $\dim\ger s\cap\ger g_1(\lambda)=4$.
	
	This phenomenon, of course, cannot occur if the symmetric superpair $(\ger g,\ger k)$ is of group type. This reflects the fact that the conditions characterising the invariant algebra may be different in the general case than one might expect from the knowledge of the group case (\emph{i.e.}~the theorems of Sergeev and Kac, Gorelik).
\end{Rem}

\subsection{The twisted action on the supersymmetric algebra}

\begin{Par}
	Let $V=V_0\oplus V_1$ be a finite-dimensional super-vector space over $\cplxs$. We define the supersymmetric algebra  $S(V)=S(V_0)\otimes\bigwedge(V_1)$. It is $\ints$-graded by total degree, as follows: $S^{k,\mathrm{tot}}(V)=\bigoplus_{p+q=k}S^p(V_0)\otimes\bigwedge^q(V_1)$. This grading is not compatible with the $\ints_2$-grading, but will of be of use to us nonetheless. 
	
	Let $U$ be another finite-dimensional super-vector space, and moreover, let $b:U\times V\to\cplxs$ be a bilinear form. Then $b$ extends to a bilinear form $S(U)\times S(V)\to\cplxs$: It is defined on linear generators by 
	\[
		b(x_1\dotsm x_m,y_1\dotsm y_n)=\delta_{mn}\cdot\sum\nolimits_{\sigma\in\ger S_n}\alpha^\sigma_{x_1,\dotsc,x_n}\cdot b(x_{\sigma(1)},y_1)\dotsm b(x_{\sigma(n)},y_n)	
	\]
	\fa $x_1,\dotsc,x_m\in U$, $y_1,\dotsc,y_n\in V$ where $\alpha=\alpha^\sigma_{x_1,\dotsc,x_n}=\pm1$ is determined by the requirement that $\alpha\cdot x_{\sigma(1)}\dotsm x_{\sigma(n)}=x_1\dotsm x_n$ in $S(V)$. If $b$ is even (resp.~odd, resp.~non-degenerate), then so is its extension. Here, recall that a bilinear form has degree $i$ if $b(V_j,V_k)=0$ whenever $i+j+k\equiv1\ (2)$.
	
	In particular, the natural pairing of $V$ and $V^*$ extends to a non-de\-ge\-ne\-rate even pairing $\Cdual0\cdot\cdot$ of $S(V)$ and $S(V^*)$. By this token, $S(V)$ embeds injectively as a subsuperspace in $\widehat S(V)=S(V^*)^*$. Its image coincides with the graded dual $S(V^*)^{*\mathrm{gr}}$ whose elements are the linear forms vanishing on $S^{k,\mathrm{tot}}(V^*)$ for $k\gg1$. 
	
	We define a superalgebra homomorphism $\partial:S(V)\to\End0{\widehat S(V^*)}$ by 
	\[
		\Cdual0p{\partial(q)\pi}=\Cdual0{pq}\pi\mathfa p,q\in S(V)\,,\,\pi\in S(V)^*
	\]
	where $\widehat S(V^*)=S(V)^*$. Clearly, $\partial(q)$ leaves $S(V^*)$ invariant. 
\end{Par}

\begin{Par}
	If $U$ is an even finite-dimensional vector space over $\cplxs$, then we have the well-known isomorphism $S(U^*)\cong\cplxs[U]$ as algebras, where $\cplxs[U]$ is the set of polynomial mappings $U\to\cplxs$. We recall that the isomorphism can be written down as follows. 
	
	The pairing $\Cdual0\cdot\cdot$ of $S(U)$ and $S(U^*)$ extends to $\widehat S(U)\times S(U^*)$. For any $d\in S(U)$, the exponential $e^d=\sum_{n=0}^\infty\frac{d^n}{n!}$ makes sense as an element of the algebra $\widehat S(U)=\prod_{n=0}^\infty S^n(U)$. Now, define a map $S(U^*)\to\cplxs[U]:p\mapsto P$ by 
	\[
		P(z)=\Cdual0{e^z}p=\textstyle\sum_{n=0}^\infty\frac1{n!}\Cdual0{z^n}{p}=\sum_{n=0}^\infty\frac1{n!}\Cdual01{\partial(z)^np}\ .
	\]
	Observe 
	\[
		\tfrac d{dt}P(z_0+tz)\big|_{t=0}=\tfrac d{dt}\Cdual0{e^{tz}e^{z_0}}p\big|_{t=0}=\Cdual0{ze^{z_0}}p\ .
	\]
	Iterating this formula, we obtain $\Cdual0{z_1\dotsm z_n}p$ for any $z_j\in U$ as a repeated directional derivative of $P$, and the map is injective. Since it preserves the grading by total degree, it is bijective because of identities of dimension in every degree. 
\end{Par}

\begin{Par}
	Let $V=V_0\oplus V_1$ be a finite-dimensional super-vector space. We apply the above to define an isomorphism $\phi:S(V^*)\to\Hom[_{S(V_0)}]0{S(V),\cplxs[V_0]}$. Here, $S(V_0)$ acts on $S(V)$ by left multiplication, and it acts on $\cplxs[V_0]$ by natural extension of the action of $V_0$ by directional derivatives:
	\[
		(\partial_zP)(z_0)=\tfrac d{dt}P(z_0+tz)\big|_{t=0}\mathfa P\in\cplxs[V_0]\,,\,z,z_0\in V_0\ .
	\]
	The isomorphism $\phi$ is given by the following prescription for $P=\phi(p)$:
	\[
		P(d;z)=(-1)^{\Abs0d\Abs0p}\Cdual0{e^z}{\partial(d)p}\mathfa p\in S(V^*)\,,\,z\in V_0\,,\,d\in S(V)\ .
	\]
	Here, note that $\widehat S(V_0)\subset\widehat S(V)$ since $S(V^*_0)$ is a direct summand of $S(V^*)$, $S(V^*)=S(V_0^*)\oplus S(V_0^*)\otimes\bigwedge^+(V_1^*)$, where $\bigwedge^+=\bigoplus_{k\sge1}\bigwedge^k$. Hence, $e^z$ may be considered as an element of $\widehat S(V)$. 	
	
	The map $\phi$ is an isomorphism as the composition of the isomorphisms
	\begin{align*}
		\Hom[_{S(V_0)}]0{S(V),\cplxs[V_0]}&\cong\Hom[_{S(V_0)}]0{S(V_0)\otimes\textstyle\bigwedge V_1,S(V_0^*)}\\
		&\cong\textstyle S(V_0^*)\otimes\bigwedge V_1^*\cong S(V^*)\ .
	\end{align*}
\end{Par}

\begin{Def}[restrhom]
	Let $(\ger g,\ger k)$ be a reductive symmetric superpair of even type, and $\ger a\subset\ger p$ an even Cartan subspace. We apply the isomorphism $\phi$ for $V=\ger p$ to define natural \emph{restriction homomorphisms}
	\[
	S(\ger p^*)\to S(\ger p_0^*):p\mapsto\bar p\nd S(\ger p^*)\to S(\ger a^*):p\mapsto\bar p\ .
	\]
	Here, $\bar p\in S(\ger p_0^*)$ (resp.~$\bar p\in S(\ger a^*)$) is defined via its associated polynomial $\bar P\in\cplxs[\ger p_0]$ (resp.~$\bar P\in\cplxs[\ger a]$) where 
	\[
		\bar P(z)=P(1;z)\nd P=\phi(p)\ .
	\]
	This is a convention we will adhere to in all that follows. 
	
	Since $\ger p_0$ is complemented by $\ger p_1$ in $\ger p$, and $\ger a$ is complemented in $\ger p_0$ by $\bigoplus_{\lambda\in\Sigma_0}\ger p_0^\lambda$, we will in the sequel consider $\ger p_0^*\subset\ger p^*$ and $\ger a^*\subset\ger p_0^*$.
\end{Def}

\begin{Par}
	Let $K$ be a connected Lie group with Lie algebra $\ger k_0$ \scth the restricted adjoint representation $\ad:\ger k_0\to\End0{\ger g}$ lifts to a homomorphism $\Ad:K\to\GL(\ger g)$. (For instance, one might take $K$ simply connected.) Then $\ger k$ (resp.~$K$) acts on $S(\ger p)$, $S(\ger p^*)$, $\widehat S(\ger p)$, $\widehat S(\ger p^*)$ by suitable extensions of $\ad$ and $\ad^*$ (resp.~$\Ad$ and $\Ad^*$) which we denote by the same symbols. Here, the sign convention for $\ad^*$ is
	\[
		\Cdual0y{\ad^*(x)\eta}=\Cdual0{[y,x]}\eta=-(-1)^{\Abs0x\Abs0y}\Cdual0{\ad(x)(y)}\eta
	\]
	\fa $x,y\in\ger g$, $\eta\in\ger g^*$. 
	
	Let $z\in\ger p_0$. We have $e^z=\sum_{k=0}^\infty\frac{z^k}{k!}\in\widehat S(\ger p)$, and this element is invertible with inverse $e^{-z}$. Define 
	\[
	u_z(x)d=\ad(x)(de^z)e^{-z}\mathfa x\in\ger k\,,\,d\in\widehat S(\ger p)\ .
	\]
	Observe that  
	\[
	\ad(x)(e^z)=\textstyle\sum_{n=0}^\infty\frac1{n!}\ad(x)(z^n)=\sum_{n=1}^\infty\frac n{n!}[x,z]z^{n-1}=[x,z]e^z\ ,
	\]	
	because $z$ is even. Hence, 
	\[
		u_z(x)d=\ad(x)(de^z)e^{-z}=[x,z]d+\ad(x)(d)\ .
	\]
	In particular, $u_z(x)$ leaves $S(\ger p)\subset\widehat S(\ger p)$ invariant. 
\end{Par}

\begin{Lem}[twistedaction-cov]
	Let $z\in\ger p_0$. Then $u_z$ defines a $\ger k$-module structure on $S(\ger p)$, and \fa $x\in\ger k$, $k\in K$, we have 
	\[
		\Ad(k)\circ u_z(x)=u_{\Ad(k)(z)}(\Ad(k)(x))\circ\Ad(k)\ .
	\]
\end{Lem}

\begin{proof}
	We clearly have 
	\[
		u_z(x)u_z(y)d=(\ad(x)\ad(y)(de^z))e^{-z}\ .
	\]
	Now $u_z$ is a $\ger k$-action because $\ad$ is a homomorphism. Similarly,
	\begin{align*}
		\Ad(k)(u_z(x)d)&=\ad(\Ad(k)(x))(\Ad(k)(d)e^{\Ad(k)(z)})e^{-\Ad(k)(z)}\\
		&=u_{\Ad(k)(z)}(\Ad(k)(x))\Ad(k)(d)\ ,
	\end{align*}
	which manifestly gives the second assertion. 
\end{proof}

\begin{Par}[twisted-action-def]
	Let $u_z$ also denote the natural extension of $u_z$ to $\Uenv0{\ger k}$. Then we may define an action $\ell$ of $\Uenv0{\ger k}$ on $\Hom[_{S(\ger p_0)}]0{S(\ger p),\cplxs[\ger p_0]}$ via 
	\[
		(\ell_vP)(d;z)=(-1)^{\Abs0v\Abs0P}P(u_z(S(v))d;z)
	\]
	\fa $P\in\Hom[_{S(\ger p_0)}]0{S(\ger p),\cplxs[\ger p_0]}$, $v\in\Uenv0{\ger k}$, $d\in S(\ger p)$, $z\in\ger p_0$. Here, we denote by $S:\Uenv0{\ger g}\to\Uenv0{\ger g}$ the unique linear map \scth $S(1)=1$, $S(x)=-x$ \fa $x\in\ger g$, and $S(uv)=(-1)^{\Abs0u\Abs0v}S(v)S(u)$ \fa homogeneous $u,v\in\Uenv0{\ger g}$ (\emph{i.e.}~the principal anti-automorphism). Compare \cite{koszul-superaction} for a similar definition in the context of the action of a supergroup on its algebra of superfunctions. 
	
	We also define 
	\[
		(L_kP)(d;z)=P(\Ad(k^{-1})(d);\Ad(k^{-1})(z))
	\]
	\fa $P\in\Hom[_{S(\ger p_0)}]0{S(\ger p),\cplxs[\ger p_0]}$, $k\in K$, $d\in S(\ger p)$, $z\in\ger p_0$.
\end{Par}

\begin{Lem}[flat-radial-fmla]
	The map $\ell$ (resp.~$L$) defines on $\Hom[_{S(\ger p_0)}]0{S(\ger p),\cplxs[\ger p_0]}$ the structure of a module over $\ger k$ (resp.~$K$) making the isomorphism $\phi$ equivariant for $\ger k$ (resp.~$K$). 
\end{Lem}

\begin{proof}
	Let $P=\phi(p)$. Then
	\begin{align*}
		(\ell_xP)(d;z)&=-(-1)^{\Abs0x\Abs0p}P(u_z(x)d;z)=-(-1)^{\Abs0d\Abs0p}\Cdual1{\ad(x)(e^zd)}p\\
		&=(-1)^{\Abs0d(\Abs0x+\Abs0p)}\Cdual1{e^zd}{\ad^*(x)(p)}=\phi\Parens1{\ad^*(x)(p)}(d;z)\ .
	\end{align*}
	
	Similarly, we check that
	\begin{align*}
		(L_kP)(d;z)&=P(\Ad(k^{-1})(d);\Ad(k^{-1})(z))\\
		&=(-1)^{\Abs0d\Abs0p}\Cdual1{e^{\Ad(k^{-1})(z)}\Ad(k^{-1})(d)}p\\
		&=(-1)^{\Abs0d\Abs0p}\Cdual1{\Ad(k^{-1})(e^zd)}p=\phi\Parens1{\Ad^*(k)(p)}(z;d)\ .
	\end{align*}
	This proves our assertion.
\end{proof}

\section{Chevalley's restriction theorem}

\subsection{The map $\gamma_z$}

\begin{Par*}
	From now on, let $(\ger g,\ger k)$ be a reductive symmetric superpair of even type, and let $\ger a\subset\ger p_0$ be an even Cartan subspace. 
\end{Par*}

\begin{Def}
	An element $z\in\ger p_0$ is called \emph{oddly regular} whenever the map $\ad(z):\ger k_1\to\ger p_1$ is surjective. 	Recall that $z\in\ger p_0$ is called \emph{regular} if $\dim\ger z_{\ger k_0}(z)=\dim\ger z_{\ger k_0}(\ger a)$. We shall call $z$ \emph{super-regular} if it is both regular and oddly regular. 
		
	Fix an even Cartan subspace $\ger a$, and let $\Sigma$ be the set of (both odd and even) restricted roots. Let $\Sigma^+\subset\Sigma$ be any subset such that $\Sigma$ is the disjoint union of $\pm\Sigma^+$. Define $\Sigma_j^\pm=\Sigma_j\cap\Sigma^\pm$ for $j\in\ints/2\ints$. Let $\bar\Sigma_1$ be the set of $\lambda\in\Sigma_1$ \scth $m\lambda\not\in\Sigma_0$ for $m=1,2$. Denote $\bar\Sigma_1^+=\bar\Sigma_1\cap\Sigma^+$. Note that $\Pi_1\in S(\ger a^*)^W$ where $\Pi_1(h)=\prod_{\lambda\in\Sigma_1}\lambda(h)$, and $W$ is the Weyl group of $\Sigma_0$. 
	
	By Chevalley's restriction theorem, restriction $S(\ger p_0^*)^{\ger k_0}\to S(\ger a^*)^W$ is a bijective map. Let $\Pi_1$ also denote the unique extension to $S(\ger p^*_0)^{\ger k_0}$ of $\Pi_1$. 
\end{Def}

\begin{Rem}
	The space $\ger p_0$ contains non-semi-simple elements, and the definitions we have given above work in this generality. However, it will suffice for our purposes to consider the set of \emph{semi-simple} super-regular elements in $\ger p_0$, by the following reasoning. 
	
	First, the set of semi-simple elements in $\ger p_0$ is Zariski dense (a linear endomorphism is semi-simple if and only if its minimal polynomial has only simple zeros). Second, the set of semi-simple elements in $\ger p_0$ equals $\Ad(K)(\ger a)$ \cite[Chapter III, Proposition~4.16]{helgason-geoman}. Thus, given any \emph{semi-simple} $z\in\ger p_0$, $z$ is oddly regular (super-regular) if and only if $\lambda(\Ad(k)(z))\neq0$ \fa $\lambda\in\Sigma_1$ ($\lambda\in\Sigma$), and for some (any) $k\in K$ \scth $\Ad(k)(z)\in\ger a$. In particular, the set of super-regular elements of $\ger a$ is the complement of a finite union of hyperplanes. Hence, the set of semi-simple super-regular elements of $\ger p_0$ is non-void and therefore Zariski dense; in particular, this holds for the set of semi-simple oddly regular elements. 
\end{Rem}

\begin{Lem}[centnondegen]
	If $z\in\ger p_0$ is semi-simple, then $\ger k_i=\ger z_{\ger k_i}(z)\oplus[z,\ger p_i]$, and the subspaces $\ger z_{\ger k_i}(z)$ and $[z,\ger p_i]$ are $b$-non-degenerate. 
\end{Lem}

\begin{proof}
	Since $\ad z$ is a semi-simple endomorphism of $\ger g$ ($\ger g$ is a semi-simple $\ger g_0$-module and $z$ is semi-simple), we have $\ger g_i=\ger z_{\ger g_i}(z)\oplus[z,\ger g_i]$. Taking $\theta$-fixed parts, we deduce $\ger k_i=\ger z_{\ger k_i}(z)\oplus\ger[z,\ger p_i]$. The summands, being $b$-orthogonal, are non-degenerate. 
\end{proof}

\begin{Par}
	Let $z\in\ger p_0$ be semi-simple and oddly regular. Let $\beta:S(\ger g)\to\Uenv0{\ger g}$ be the supersymmetrisation map. Let
	\[
		\Gamma_z:\textstyle\bigwedge(\ger p_1)\otimes S(\ger p_0)\to S(\ger p):q\otimes p\mapsto u_z\Parens1{\beta([z,q])}p
	\]
	on elementary tensors and extend linearly. 
\end{Par}

\begin{Prop}[radial-part]
Let $z$ be oddly regular and semi-simple. Then $\Gamma_z$ is bijective, and $\gamma_z=(\eps\otimes 1)\circ\Gamma_z^{-1}:S(\ger p)\to S(\ger p_0)$ satisfies 
\[
\gamma_{\Ad(k)(z)}\circ\Ad(k)=\Ad(k)\circ\gamma_z\mathfa k\in K\ .
\] 
Here $\eps:\bigwedge(\ger p_1)\to\cplxs$ is the unique unital algebra homomorphism.

Moreover, on $S^{m,\mathrm{tot}}(\ger p)$, $\Pi_1(z)^m\gamma_z$ is polynomial in $z$, \emph{i.e.}~it extends to an element $\Pi_1(\cdot)^m\gamma_\cdot$ of the space $\cplxs[\ger p_0]\otimes\Hom0{S^{m,\mathrm{tot}}(\ger p),S(\ger p_0)}$. 
\end{Prop}

\begin{proof}
	By the assumption on $z$, $\ad z:\ger p_1\to[z,\ger p_1]$ is bijective. Moreover, $\Gamma_z$ respects the filtrations by total degree, and the degrees of these filtrations are equidimensional by the assumption. Hence, $\Gamma_z$ will be bijective once it is surjective. In degree zero, $\Gamma_z$ is the identity. We proceed to prove the surjectivity in higher degrees by induction.

	By assumption, $\ad z:[z,\ger p_1]\to\ger p_1$ is also bijective (since its kernel is $\ger z_{\ger k_1}(z)\cap[z,\ger p_1]$, which is $0$ by \thmref{Lem}{centnondegen}). Let $y_1,\dotsc,y_m\in\ger p_1$, $y_1',\dotsc,y_n'\in\ger p_0$. Let $x_j\in\ger p_1$ \scth $[[z,x_j],z]=y_j$. We find
	\[
		\Gamma_z(x_1\dotsm x_m\otimes y_1'\dotsm y_n')\equiv y_1\dotsm y_my_1'\dotsm y_n'\quad\Parens1{\textstyle\bigoplus\nolimits_{k<m+n}S^{k,\mathrm{tot}}(\ger p)}\ ,
	\]
	so the first assertion follows by induction. 
	
	As to the covariance property, observe first that we have the identity $\Ad(k)([z,\ger p_1])=[\Ad(k)(z),\ger p_1]$. Moreover, 
	\begin{align*}
		(\Ad(k)\circ\gamma_z)(\Gamma_z(v\otimes d))&=\eps(v)\Ad(k)(d)=\eps(\Ad(k)(v))\Ad(k)(d)\\
		&=\gamma_{\Ad(k)(z)}\Parens1{\Gamma_{\Ad(k)(z)}(\Ad(k)(v)\otimes\Ad(k)(d))}\\
		&=\gamma_{\Ad(k)(z)}\Parens1{u_{\Ad(k)(z)}(\Ad(k)(\beta([z,v])))\Ad(k)(d)}\\
		&=\gamma_{\Ad(k)(z)}\Parens1{\Ad(k)(u_z(\beta([z,v]))(d))}\\
		&=(\gamma_{\Ad(k)(z)}\circ\Ad(k))(\Gamma_z(v\otimes d))
	\end{align*}
	\fa $v\in\bigwedge(\ger p_1)$ and $d\in S(\ger p_0)$, by \thmref{Lem}{twistedaction-cov}.
	
	To show that $\Pi_1(z)^m\gamma_z:S^{m,\mathrm{tot}}(\ger p)\to S(\ger p_0)$ is given by the restriction of a polynomial function, we remark that its domain of definition---the set $U$ of semi-simple oddly regular elements in $\ger p_0$---is Zariski dense. We need only prove that $f:U\to\Hom0{\ger p_1,\ger k_1}$, $f(z)=\Pi_1(z)(\ad z)^{-1}$, is polynomial in $z$, where we consider $(\ad z)^{-1}:\ger p_1\to[z,\ger p_1]$ as a linear map $\ger p_1\to\ger k_1$. 
		
	Thus, let $z\in\ger p_0$ be semi-simple and oddly regular. It is contained in some even Cartan subspace $\ger a$ (say). We have $\ger z_{\ger k_1}(\ger a)=\ger m_1=\ger k_1\cap[z,\ger p_1]^\perp$ by \thmref{Lem}{centnondegen} and $(\ger k_1\cap\ger m_1^\perp)\oplus\ger p_1=\bigoplus_{\lambda\in\Sigma_1^+}\ger g_{1,\ger a}^\lambda$. If $x=u+v\in\ger g^\lambda_{1,\ger a}$, and $u\in\ger k_1$, $v\in\ger p_1$, then $[z,u]=\lambda(z)v$. It follows that $\Pi_1(z)(\ad z)^{-1}$ depends polynomially on $z$, proving our claim. 
\end{proof}

\begin{Prop}[flat-radial]
	Let $p\in S(\ger p^*)^{\ger k}$. Then $P(d;z)=P(\gamma_z(d);z)$ \fa oddly regular and semi-simple $z\in\ger p_0$ and $d\in S(\ger p)$. 
\end{Prop}

\begin{proof}
	Fix an oddly regular $z\in\ger p_0$, and let $x_1,\dotsc,x_n\in\ger p_1$. By \thmref{Lem}{flat-radial-fmla}, we find for $n>0$ 
	\[
		P\Parens1{\Gamma_z(x_1\dotsm x_n\otimes q);z}=(\ell_{\beta([z,x_1\dotsm x_n])}P)(q;z)=0\ .
	\]
	Since $d-\gamma_z(d)\in\Gamma_z(\bigwedge^+(\ger p_1)\otimes S(\ger p_0))$, where $\bigwedge^+(\ger p_1)$ denotes the kernel of $\eps:\bigwedge(\ger p_1)\to\cplxs$ (\emph{i.e.}, the set of elements without constant term), the assertion follows immediately.  
\end{proof}

\begin{Cor}[res-inj]
	Let $(\ger g,\ger k)$ be a reductive symmetric superpair of even type. The algebra homomorphism $p\mapsto\bar p:I(\ger p^*)=S(\ger p^*)^{\ger k}\to S(\ger p_0^*)$ is injective. In particular, $I(\ger p^*)$ is commutative and purely even. 
\end{Cor}

\begin{proof}
	Let $p\in I(\ger p^*)$. Assume that $\bar p=0$. Let $d\in S(\ger p)$. For all $z\in\ger p_0$ which are oddly regular and semi-simple, 
	\[
		P(d;z)=P(\gamma_z(d);z)=[\partial_{\gamma_z(d)}\bar P](z)=0\ ,
	\]
	by \thmref{Prop}{flat-radial}. It follows that $P(d;-)=0$ on $\ger p_0$, since it is a polynomial. Since $d$ was arbitrary, we have established our contention. 
\end{proof}

\begin{Rem}
	The statement of the Corollary can, of course, be deduced by applying the inverse function theorem for supermanifolds, as in \cite[Proposition 1.1]{sergeev-invpol}. Nonetheless, we find it instructive to give the above proof based on the map $\gamma_z$, as it illustrates the approach we will take to determine the image of the restriction map. 
\end{Rem}

\subsection{Proof of Theorem (A)} 

\begin{Par}
	Let $(\ger g,\ger k)$ be a reductive symmetric superpair of even type, and let $\ger a$ be an even Cartan subspace. We denote by $\ger a'$ the set of super-regular elements of $\ger a$. Let $\mathcal R$ be the algebra of differential operators on $\ger a$ with rational coefficients which are non-singular on $\ger a'$. For any $z\in\ger a'$ and any $D\in\mathcal R$, let $D(z)$ be the local expression of $D$ at $z$. This is defined by the requirement that $D(z)$ be a differential operator with constant coefficients, and 
	\[
		(Df)(z)=(D(z)f)(z)\mathfa z\in\ger a'\ ,
	\]
	and all regular functions $f$. 
	
	We associate to $\Sigma\subset\ger a^*$, the restricted root system of $\ger g:\ger a$, the subset $\mathcal R_\Sigma=\bigcup_{\lambda\in\bar\Sigma_1^+}\mathcal R_\lambda\subset\mathcal R$ where
	\[
		\mathcal R_\lambda=\Set1{D\in\mathcal R}{\exists\,d\in S(\ger p_1^\lambda)\colon D(z)=\gamma_z(d)\text{ \fa}z\in\ger a'}\ .
	\]
	\emph{I.e.}, $\mathcal R_\Sigma$ consists of those differential operators which are given as radial parts of operators with constant coefficients on the $\ger p$-projections $\ger p_1^\lambda$ of the restricted root spaces for the $\lambda\in\bar\Sigma_1^+$. For any $D\in\mathcal R$, let the \emph{domain} $\dom D$ be the set of all $p\in\cplxs[\ger a]$ \scth $Dp\in\cplxs[\ger a]$. 
	
	As we shall see, the image of the restriction map is the set of $W$-invariant polynomials in the common domain of $\mathcal R_\Sigma$. We will subsequently determine $\mathcal R_\Sigma$ in order to describe this common domain in more explicit terms. 
\end{Par}

\begin{Th}[super-chevalley]
	The restriction homomorphism $I(\ger p^*)\to S(\ger a^*)$ from \thmref{Def}{restrhom} is a bijection onto the subspace $I(\ger a^*)=S(\ger a^*)^W\cap\bigcap_{D\in\mathcal R_\Sigma}\dom D$. 
\end{Th}

\begin{Par*}
	The \emph{proof} of the Theorem requires a little preparation. 
\end{Par*}

\begin{Lem}[doubleprojection-trick]
	Let $q\in S(\ger p_0^*)^K$, $Q=\phi(q)$, and $z\in\ger p_0$ be super-regular and semi-simple. \Fa $x\in\ger k$, and $w\in S(\ger p)$, we have
	\[
		Q\Parens1{\gamma_z(u_z(x)w);z}=0\ .
	\]
\end{Lem}

\begin{proof}
	There is no restriction to generality in supposing $z\in\ger a'$, so that $\ger z_{\ger k}(z)=\ger z_{\ger k}(\ger a)=\ger m$ and $\ger z_{\ger k_0}(z)=\ger z_{\ger k_0}(\ger a)=\ger m_0$. We define linear maps
	\[
		\gamma_z':S(\ger p_0)\to S(\ger a)\nd\gamma_z'':S(\ger p)\to S(\ger a)
	\]
	by the requirements that $v-\gamma_z'(v)\in u_z(\ger m_0^\perp\cap\ger k_0)(S(\ger p_0))$ \fa $v\in S(\ger p_0)$ and $w-\gamma_z''(w)\in u_z(\ger m^\perp\cap\ger k)(S(\ger p))$ \fa $w\in S(\ger p)$. (That such maps exist and are uniquely defined by these properties follows in exactly the same way as for \thmref{Prop}{radial-part}. We remark that $[z,\ger p_i]=\ger k_i\cap\ger m_i^\perp$ by \thmref{Lem}{centnondegen}.) Then 
	\begin{align*}
		w&-\gamma_z'(\gamma_z(w))=w-\gamma_z(w)+\gamma_z(w)-\gamma_z'(\gamma_z(w))\\
		&\in u_z(\ger m_1^\perp\cap\ger k_1)(S(\ger p))+u_z(\ger m_0^\perp\cap\ger k_0)(S(\ger p_0))\subset u_z(\ger m^\perp\cap\ger k)(S(\ger p))
	\end{align*}
	\fa $w\in S(\ger p)$, where $\ger m_1=\ger z_{\ger k_1}(\ger a)$. This shows that $\gamma_z''=\gamma_z'\circ\gamma_z$. 
	
	Moreover, by the $K$-invariance of $q$, we have $Q(v;z)=Q(\gamma_z'(v);z)$ \fa $v\in S(\ger p_0)$. We infer  
	\[
		Q\Parens1{\gamma_z(u_z(x)w);z}=Q\Parens1{\gamma_z''(u_z(x)w);z}=0\mathfa x\in\ger m^\perp\cap\ger k\ ,\ w\in S(\ger p)
	\]
	since $u_z(x)w\in u_z(\ger m^\perp\cap\ger k)(S(\ger p))$ belongs to $\ker\gamma_z''$. 
	
	Next, we need to consider the case of $x\in\ger m$. Then $\ad(x):S(\ger p)\to S(\ger p)$ annihilates the subspace $S(\ger a)$, and moreover, $\ad(x)(e^z)=0$. From this we find \fa $y\in\ger m^\perp\cap\ger k$, $d\in S(\ger p)$
	\begin{align*}
		\ad(x)\Parens1{u_z(y)(d)}&=(\ad(x)\ad(y)(de^z))e^{-z}\\
		&=(\ad([x,y])(de^z))e^{-z}+(-1)^{\Abs0x\Abs0y}\ad(y)(\ad(x)(d)e^z)e^{-z}\\
		&=u_z([x,y])d+(-1)^{\Abs0x\Abs0y}u_z(y)\ad(x)(d)\ .
	\end{align*}
	Since $\ger m$ is a subalgebra and $b$ is $\ger k$-invariant, $\ger m^\perp\cap\ger k$ is $\ger m$-invariant. Hence, the above formula shows that $\ker\gamma_z''=u_z(\ger m^\perp\cap k)(S(\ger p))$ is $\ad(x)$-invariant.
	
	By the definition of $\gamma_z''$, we find that 
	\[
		\gamma_z''(\ad(x)d)=\ad(x)\gamma_z''(d)=0\mathfa x\in\ger m\,,\,d\in S(\ger p)\ .
	\]
	Reasoning as above, we see that 
	\[
	Q(\gamma_z(u_z(x)d);z)=Q(\gamma_z(\ad(x)d);z)=0\mathfa x\in\ger m\,,\,d\in S(\ger  p)\ .
	\]
	Since $\ger k=\ger m\oplus(\ger m^\perp\cap\ger k)$, this proves the lemma. 
\end{proof}

\begin{Par*}
	Let $\ger p_0'$ be the set of semi-simple super-regular elements in $\ger p_0$. Recall the polynomial $\Pi_1$, and consider the localisation $\cplxs[\ger p_0]_{\Pi_1}$. Let $q\in S(\ger p_0^*)^K$, $Q=\phi(q)$, and define 
	\[
		P(v;z)=Q(\gamma_z(v);z)\mathfa v\in S(\ger p)\ ,\ z\in\ger p_0'\ .
	\]
	By \thmref{Prop}{radial-part}, $P\in\Hom0{S(\ger p),\cplxs[\ger p_0]_{\Pi_1}}$. We remark that the $\ger k$-action $\ell$ defined in \thmref{Par}{twisted-action-def} extends to $\Hom0{S(\ger p),\cplxs[\ger p_0]_{\Pi_1}}$, by the same formula. 
\end{Par*}

\begin{Lem}[sp0-linear]
	Retain the above assumptions. Then $P$ is $S(\ger p_0)$-linear and $\ger k$-invariant, \emph{i.e.}~$P\in\Hom[_{S(\ger p_0)}]0{S(\ger p),\cplxs[\ger p_0]_{\Pi_1}}^{\ger k}$.
\end{Lem}

\begin{proof}
	By \thmref{Lem}{doubleprojection-trick}, $P$ is $\ger k$-invariant. It remains to prove that $P$ is $S(\ger p_0)$-linear. To that end, we first establish that $P$ is $K$-equivariant as linear map $S(\ger p)\to\cplxs[\ger p_0]_{\Pi_1}$. Since $q$ is $K$-invariant, 
	\begin{align*}
		P\Parens1{\Ad(k)(v);\Ad(k)(z)}&=Q\Parens1{\gamma_{\Ad(k)(z)}(\Ad(k)(v));\Ad(k)(z)}\\
		&=Q\Parens1{\Ad(k)(\gamma_z(v));\Ad(k)(z)}\\
		&=Q(\gamma_z(v);z)=P(v;z)\ .
	\end{align*}

	Next, fix $z\in\ger p_0'$. Then $S(\ger p)=S(\ger p_0)\oplus u_z(\ger z_{\ger k_1}(z)^\perp\cap\ger k_1)(S(\ger p))$ where the second summand equals $\ker\gamma_z$. We may check the $S(\ger p_0)$-linearity on each summand separately. 

	For $v\in S(\ger p_0)$, we have $P(v;z)=Q(v;z)$, so for any $y\in\ger p_0$ 
	\[
		[\partial_yP(v;-)](z)=[\partial_yQ(v;-)](z)=Q(yv;z)=P(yv;z)\ .
	\]
	
	We are reduced to considering $v=u_z(x)v'$ where $x\in\ger z_{\ger k_1}(z)^\perp\cap\ger k_1$ and $v'\in S(\ger p)$. We may assume w.l.o.g.~$z\in\ger a$ (since $z$ is semi-simple), so that $\ger z_{\ger k_1}(z)=\ger z_{\ger k_1}(\ger a)=\ger m_1$. By our assumption on $z$, $\ger p_0=\ger a\oplus[\ger k_0,z]$, and we may consider $y$ in each of the two summands separately. 
	
	Let $y\in\ger a$. For sufficiently small $t$, we have $z+ty\in\ger a'=\ger a\cap\ger p_0'$, so that $\ger z_{\ger k_1}(z+ty)=\ger m_1=\ger z_{\ger k_1}(z)$. Hence, $\gamma_{z+ty}(u_{z+ty}(x)v')=0$. By the chain rule,  
	\[
		0=\tfrac d{dt}\gamma_{z+ty}(u_{z+ty}(x)v')\big|_{t=0}=d\gamma_\cdot(v)_z(y)+\gamma_z\Parens1{\tfrac d{dt}u_{z+ty}(x)v'\big|_{t=0}}\ ,
	\]
	Since $\tfrac d{dt}u_{z+ty}(x)v'\big|_{t=0}=[x,y]v'$, we have 
	\[
		d\gamma_\cdot(v)_z(y)=-\gamma_z(\tfrac d{dt}u_{z+ty}(x)v'\big|_{t=0})=\gamma_z([y,x]v')\ .
	\]
	
	Moreover, as operators on $S(\ger p)$, 
	\[
		[y,u_z(x)]=y[x,z]+y\ad(x)-[x,z]y-\ad(x)y=[y,x]\ ,
	\]
	and thus $yv=yu_z(x)v'\equiv[y,x]v'$ modulo $\ker\gamma_z$. We conclude  
	\[
		d\gamma_\cdot(v)_z(y)=\gamma_z([y,x]v')=\gamma_z(yv)=\gamma_z(yv)-y\gamma_z(v)
	\]
	since $\gamma_z(v)=0$. 	Hence, 
	\[
		[\partial_yP(v;-)](z)=Q\Parens1{d\gamma_\cdot(v)_z(y)+y\gamma_z(v);z}=Q\Parens1{\gamma_z(yv);z}=P(yv;z)\ .
	\]

	Now let $y=[u,z]$ where $u\in\ger k_0$. We may assume that $u\perp\ger z_{\ger k_0}(z)$. Define $k_t=\exp tu$. Then by the $K$-invariance of $P$, 
	\begin{align*}
		[\partial_yP(v;-)](z)&=\tfrac d{dt}P\Parens1{v;\Ad(k_t)(z)}\big|_{t=0}
		=\tfrac d{dt}P\Parens1{\Ad(k_t^{-1})(v);z}\big|_{t=0}\\
		&=-P\Parens1{\ad(u)(v);z}=P(yv;z)-P\Parens1{u_z(u)v;z}=P(yv;z)
	\end{align*}
	where in the last step, we have used \thmref{Lem}{doubleprojection-trick}. 
\end{proof}

\begin{proof}[\protect{Proof of \thmref{Th}{super-chevalley}}]
	The restriction map is injective by \thmref{Cor}{res-inj} and Chevalley's restriction theorem for $\ger g_0$. By the latter, the image lies in the set of $W$-invariants. Let $\bar p\in S(\ger a^*)$ be the restriction of $p\in I(\ger p^*)$, and $P=\phi(p)$. For any $d\in S(\ger p)$, and $D\in\mathcal R_\Sigma$ given by $D(z)=\gamma_z(d)$, we have by \thmref{Prop}{flat-radial}
	\[
		(D\bar p)(z)=(\partial_{\gamma_z(d)}\bar P)(z)=P(\gamma_z(d);z)=P(d;z)\mathfa z\in\ger a'\ .
	\]
	The result is clearly polynomial in $z$, so $\bar p\in\dom D$. This shows that the image of the restriction map lies in $I(\ger a^*)$. 
	
	Let $r\in I(\ger a^*)$. By Chevalley's restriction theorem, there exists a unique $q\in I(\ger p_0^*)=S(\ger p_0^*)^K$ \scth $Q(h)=R(h)$ \fa $h\in\ger a$. 
	
	Next, recall that for $d\in S(\ger p)$ and $z\in\ger p_0'$: 
	\[
		P(d;z)=Q(\gamma_z(d);z)\ .
	\]
	By \thmref{Lem}{sp0-linear}, $P\in\Hom[_{S(\ger p_0)}]0{S(\ger p),\cplxs[\ger p_0]_{\Pi_1}}^{\ger k}$. Hence, $P$ will define an element $p\in I(\ger p^*)$ by virtue of the isomorphism $\phi$, as soon as it is clear that, as a linear map $S(\ger p)\to\cplxs[\ger p_0]_{\Pi_1}$, it takes its values in $\cplxs[\ger p_0]$. 
	
	We only have to consider $z$ in the Zariski dense set $\ger p_0'$. The function $\Pi_1(z)^k\cdot P(d;z)$ depends polynomially on $z$, where we assume $d\in S^{\sle k,\mathrm{tot}}(\ger p)$.  To prove that $P$ has polynomial values, it will suffice (by the removable singularity theorem and the conjugacy of Cartan subspaces) to prove that $P(d;h)$ is bounded as $h\in\ger a'=\ger a\cap\ger p_0'$ approaches one of the hyperplanes $\lambda^{-1}(0)$ where $\lambda\in\Sigma_1^+$ is arbitrary. Since $r$ is $W$-invariant, $r-r_0$ (where $r_0$ is the constant term of $r$) vanishes on $\lambda^{-1}(0)$ if a multiple of $\lambda$ belongs to $\Sigma_0^+$. Such a multiple could only be $\pm\lambda,\pm2\lambda$. Hence, it will suffice to consider $\lambda\in\bar\Sigma_1^+$. By definition, $2\lambda\not\in\Sigma$. 
	
	Consider $P(d;h)$ as a map linear in $d$, and let $N_h=\ker P(-;h)$. Let $d\in S^{\sle k,\mathrm{tot}}(\ger p)$. Assume that $d=zd'$ where $z$ is defined by $x=y+z$, $y\in\ger k$, $z\in\ger p$, \fs $x\in\ger g^\mu_{\ger a}$ and $\mu\in\Sigma^+$, $\mu\neq\lambda$. Then, modulo $N_h$, 
	\[
		d=zd'\equiv zd'+\frac{u_h(y)d'}{\mu(h)}=zd'+\frac{[y,h]d'}{\mu(h)}+\frac{\ad(y)(d')}{\mu(h)}=\frac{\ad(y)(d')}{\mu(h)}\ .
	\]
	The root $\mu$ is not proportional to $\lambda$ and the total degree of $\ad(y)(d')$ is strictly less than that of $d$. By induction, modulo $N_h$, 
	\[
		d\equiv\frac{\tilde d}{\Pi_{\mu\in\Sigma^+\setminus\lambda}\,\mu(h)^k}
	\]
	\fs $\tilde d$ which lies in the subalgebra of $S(\ger p)$ generated by $\ger a\oplus\ger p_1^\lambda$, and depends polynomially on $h$ and linearly on $d\in S^{\sle k,\mathrm{tot}}(\ger p)$. 
	
	Hence, the problem of showing that $P(d;h)$ remains bounded as $h$ approaches $\lambda^{-1}(0)$ is reduced to the case of $d\in S(\ger a\oplus\ger p_1^\lambda)$. For $d\in S(\ger p_1^\lambda)$, the polynomiality of $P(d;-)$ immediately follows from the assumption on $r$. 
If $d=d'd''$ where $d'\in S(\ger a)$ and $d''\in S(\ger p_1^\lambda)$, then $P(d;z)=[\partial(d')P(d'';-)](z)$ since $P$ is $S(\ger p_0)$-linear. But $P(d'';-)\in\cplxs[\ger p_0]$ and this space is $S(\ger p_0)$-invariant, so $P(d;-)\in\cplxs[\ger p_0]$. 

	Therefore, there exists $p\in I(\ger p^*)$ \scth $P=\phi(p)$. By its definition, it is clear that $p$ restricts to $r$, so we have proved the theorem.
\end{proof}

\subsection{Proof of Theorem (B)}

\begin{Par}[symplbasis]
	In order to give a complete description of the image of the restriction map, we need to compute the radial parts $\gamma_h(d)$ for $d\in S(\ger p_1^\lambda)$ and $h\in\ger a'$ explicitly. First, let us choose bases of the spaces $S(\ger p_1^\lambda)$. 

	Let $\lambda\in\Sigma^+_1$. By \thmref{Prop}{res-superroot}~(v) we may choose $b^\theta$-symplectic bases $y_i,\tilde y_i\in\ger k_1^\lambda$, $z_i,\tilde z_i\in\ger p_1^\lambda$, $i=1,\dotsc,\frac12m_{1,\lambda}$, $m_{1,\lambda}=\dim\ger g_{1,\ger a}^\lambda$. \emph{I.e.},
	\[
		b(y_i,\tilde y_j)=b(\tilde z_j,z_i)=\delta_{ij}\,,\,b(y_i,y_j)=b(\tilde y_i,\tilde y_j)=b(z_i,z_j)=b(\tilde z_i,\tilde z_j)=0\ .
	\]
	We may impose the conditions $x_i=y_i+z_i,\tilde x_i=\tilde y_i+\tilde z_i\in\ger g_{1,\ger a}^\lambda$, so that  
	\[	
		[h,y_i]=\lambda(h)z_i\,,\,[h,\tilde y_i]=\lambda(h)\tilde z_i\,,\,[h,z_i]=\lambda(h)y_i\,,\,[h,\tilde z_i]=\lambda(h)\tilde y_i
	\]
	\fa $h\in\ger a$. (Compare \thmref{Prop}{res-superroot}~(iv).) 
	
	Given partitions $I=(i_1<\dotsm<i_k)$, $J=(j_1<\dotsm<j_\ell)$, we define monomials $z_I\tilde z_J=z_{i_1}\dotsm z_{i_k}\tilde z_{j_1}\dotsm\tilde z_{j_\ell}$ in $S(\ger p_1^\lambda)=\bigwedge(\ger p_1^\lambda)$. They form a basis of $S(\ger p_1^\lambda)$. 
\end{Par}

\begin{Lem}[radial-rank1-odd]
	Fix $\lambda\in\bar\Sigma_1^+$. Let $h\in\ger a$ be oddly regular, $I,J$ be multi-indices where $k=\Abs0I$, $\ell=\Abs0J$, and let $m$ be a non-negative integer. Modulo $\ker\gamma_h$,
	\[
		z_I\tilde z_JA_\lambda^m\equiv\begin{cases}0&I\neq J\ ,\\
		A_\lambda^m & I=J=\vvoid\ ,\\
		(-1)^kz_{I'}\tilde z_{I'}\textstyle\sum_{j=0}^m(-1)^j\tfrac{\lambda(A_\lambda)^j}{\lambda(h)^{j+1}}(m)_jA_\lambda^{m+1-j}&I=J=(i<I')\ ,\end{cases}
	\]
	where $(m)_j$ is the falling factorial $m(m-1)\dotsm(m-j+1)$, and $(m)_0=1$.
	\end{Lem}

\begin{proof}
	 For $k=\ell=0$, there is nothing to prove. We assume that $k>0$ or $\ell>0$, and write $I=(i<I')$ if $k>0$, $J=(j<J')$ if $\ell>0$. We claim that modulo $\ker\gamma_h$,
	\[
		z_I\tilde z_JA_\lambda^m\equiv\begin{cases}0&k\neq\ell\text{ or }i\neq j\ ,\\
		(-1)^kz_{I'}\tilde z_{J'}\textstyle\sum_{n=0}^m(-1)^n\tfrac{\lambda(A_\lambda)^n}{\lambda(h)^{n+1}}(m)_nA_\lambda^{m+1-n}&i=j\ .\end{cases}
	\]
	We argue by induction on $\max(k,\ell)$. There will also be a sub-induction on the integer $m$. First, we assume that $k>0$, and compute 
	\[
		z_I\tilde z_JA_\lambda^m\equiv z_iz_{I'}\tilde z_JA_\lambda^m+\tfrac1{\lambda(h)}u_h(y_i)(z_{I'}\tilde z_JA_\lambda^m)=\tfrac1{\lambda(h)}\ad(y_i)(z_{I'}\tilde z_JA_\lambda^m)\ .
	\]
	
	For any $q$, we have 
	\[
		b\Parens1{[y_i,z_q],h'}=-\lambda(h')b(y_i,y_q)=0\mathfa h'\in\ger a\ ,
	\]
	so $b([y_i,z_q],\ger a)=0$, and $[y_i,z_q]\in\ger p_0$. Hence $[y_i,z_q]\in\ger g_{0,\ger a}^{2\lambda}\oplus\ger g_{0,\ger a}^{-2\lambda}=0$. Similarly, for $i\neq q$, we have $[y_i,\tilde z_q]=0$. Now, assume that $i\sle J$. Then 
	\begin{align*}
	z_I\tilde z_JA_\lambda^m&\equiv(-1)^{k-1}\tfrac1{\lambda(h)}z_{I'}\ad(y_i)(\tilde z_JA_\lambda^m)\\
	&=(-1)^{k-1}\tfrac1{\lambda(h)}[y_i,\tilde z_j]z_{I'}\tilde z_{J'}A_\lambda^m-m\tfrac{\lambda(A_\lambda)}{\lambda(h)}z_I\tilde z_JA_\lambda^{m-1}\tag{$*$}
	\end{align*}
	since $[y_i,A_\lambda^m]=-m\lambda(A_\lambda)z_iA_\lambda^{m-1}$. As it stands, equation ($*$) only holds for $\ell>0$, but if we take the first summand to be $0$ if $\ell=0$, then it is also true in the latter case.
	
	If $\ell>0$ and $i<J$, then the first summand also vanishes, and arguing by induction on $m$, we find 
	\[
	z_I\tilde z_JA_\lambda^m\equiv(-1)^mm!\tfrac{\lambda(A_\lambda)^m}{\lambda(h)^m}z_I\tilde z_J=(-1)^{m+k-1}m!\tfrac{\lambda(A_\lambda)^m}{\lambda(h)^{m+1}}[y_i,\tilde z_j]z_{I'}\tilde z_J=0\ .
	\]
	 Virtually the same reasoning goes through for $\ell=0$. In particular, whenever $\gamma_h(z_I\tilde z_JA_\lambda^m)\neq0$ and $k>0$, then $i\sle J$ implies $\ell>0$ and $i=j$.
		
	If $\ell>0$ and $j\sle I$, then we observe that $z_I\tilde z_J=(-1)^{k\ell}\tilde z_Jz_I$. Formally exchanging the letters $z_s$ and $\tilde z_s$ in the above equations, and reordering all terms in the appropriate fashion, we obtain
	\[\tag{$**$}
		z_I\tilde z_JA_\lambda^m\equiv(-1)^k\tfrac1{\lambda(h)}[\tilde y_j,z_i]z_{I'}\tilde z_{J'}A_\lambda^m-m\tfrac{\lambda(A_\lambda)}{\lambda(h)}z_I\tilde z_JA_\lambda^{m-1}\ ,
	 \]
	 because $k\ell+\ell-1+(k-1)(\ell-1)=k(2\ell-1)\equiv k\ (2)$. Arguing as above, the right hand side of equation ($**$) is equivalent to $0$ modulo $\ker\gamma_h$ if $k=0$ or $j<I$. Therefore, $\gamma_h(z_I\tilde z_JA_\lambda^m)$ vanishes unless $k,\ell>0$ and $i=j$. 
	 
	 We consider the case of $k,\ell>0$ and $i=j$. Since $[y_i,\tilde z_i]-[\tilde y_i,z_i]=-2A_\lambda$ by standard arguments, we find, by adding equations ($*$) and ($**$), 
	\[
		z_I\tilde z_JA_\lambda^m\equiv(-1)^k\tfrac1{\lambda(h)}z_{I'}\tilde z_{J'}A_\lambda^{m+1}-m\tfrac{\lambda(A_\lambda)}{\lambda(h)}z_I\tilde z_JA_\lambda^{m-1}\ .
	\]
	We may now apply this formula recursively to the second summand, to conclude
	\[
	z_I\tilde z_JA_\lambda^m\equiv(-1)^kz_{I'}\tilde z_{J'}\textstyle\sum_{n=0}^m(-1)^n\tfrac{\lambda(A_\lambda)^n}{\lambda(h)^{n+1}}(m)_nA_\lambda^{m+1-n}\ .
	\]

	By induction on $\max(k,\ell)$, the right hand side belongs to $\ker\gamma_h$ unless $k=\ell$. We have proved our claim, and thus, we arrive at the assertion of the lemma. 
\end{proof}

\begin{Par}	
	Fix $\lambda\in\bar\Sigma_1^+$ and $h\in\ger a'$. Let $I=(i_1<\dotsm<i_k)$ and $1\sle\ell\sle k$. Set $I'=(i_{\ell+1}<\dotsm<i_k)$. Let 
	\[
		\eps^k_\ell=(-1)^{\sum_{j=k-\ell+1}^kj}=(-1)^{\frac\ell2(2k-\ell+1)}\ .
	\]
	We claim that there are $b_{s\ell}\in\nats$, $s<\ell$, $b_{01}=1$, such that, modulo $\ker\gamma_h$, 
	\[
		z_I\tilde z_I\equiv\eps^k_\ell z_{I'}\tilde z_{I'}\textstyle\sum_{j=0}^{\ell-1} b_{j\ell}\tfrac{(-\lambda(A_\lambda))^j}{\lambda(h)^{\ell+j}}A_\lambda^{\ell-j}\ .\tag{$***$}
	\]
	
	The case $\ell=1$ has already been established. To prove the inductive step, let $I''=(i_\ell,\dotsc,i_k)=(i_\ell<I')$, and $J=(i_0<I)$. We compute
	\begin{align*}
		z_J\tilde z_J&\equiv\eps_\ell^{k+1}z_{I''}\tilde z_{I''}\textstyle\sum_{j=0}^{\ell-1}b_{j\ell}\tfrac{(-\lambda(A_\lambda))^j}{\lambda(h)^{\ell+j}}A_\lambda^{\ell-j}\\
	&\equiv(-1)^{k-\ell+1}\eps_\ell^{k+1}z_{I'}\tilde z_{I'}\textstyle\sum_{s=0}^\ell\sum_{j=0}^{\min(s,\ell-1)}(\ell-j)_{s-j}b_{j\ell}\tfrac{(-\lambda(A_\lambda))^s}{\lambda(h)^{\ell+1+s}}A_\lambda^{\ell+1-s}\ ,
	\end{align*}
	so 
	\[
		b_{s,\ell+1}=\textstyle\sum_{j=0}^{\min(s,\ell-1)}(\ell-j)_{s-j}b_{j\ell}=\frac1{(\ell-s)!}\sum_{j=0}^{\min(s,\ell-1)}(\ell-j)!b_{j\ell}\ .
	\]
	This proves our claim, where the constants $b_{s\ell}$ obey the recursion relation set out above. 
	
	To solve this recursion, we claim that 
	\[
		b_{s\ell}=\frac{(\ell-1+s)!}{2^s(\ell-1-s)!s!}\mathfa 0\sle s<\ell\ .
	\]
	This is certainly the case for $\ell=1$. By induction, \fa $0\sle s\sle\ell$, $\ell\sge1$, 
	\[
		b_{s,\ell+1}=\tfrac1{(\ell-s)!}\textstyle\sum_{j=0}^{\min(s,\ell-1)}(\ell-j)\tfrac{(\ell-1+j)!}{2^jj!}\ .
	\]
	As is easy to show by induction, $\sum_{j=0}^N(\ell-j)\tfrac{(\ell-1+j)!}{2^jj!}=\tfrac{(\ell+N)!}{2^NN!}$. Hence, 
	\[
		b_{s,\ell+1}=\begin{cases}\frac{(\ell+s)!}{2^s(\ell-s)!s!}&0\sle s<\ell\\\frac{(2\ell-1)!}{2^{\ell-1}(\ell-1)!}=\frac{(2\ell)!}{2^\ell\ell!} & s=\ell\end{cases}
	\]
	which establishes the claim. 
	
	Setting $\ell=k=\Abs0I$ in $(***)$, we obtain the following lemma.
\end{Par}

\begin{Lem}[radialfmla]
	Fix $\lambda\in\bar\Sigma_1^+$. Let $h\in\ger a$ be oddly regular, $I$ be a multi-index where $k=\Abs0I$. Then 
	\[
		\gamma_h(z_I\tilde z_I)=(-1)^{\frac{k(k+1)}2}\textstyle\sum_{j=0}^{k-1}\tfrac{(k-1+j)!}{2^j(k-1-j)!j!}\tfrac{(-\lambda(A_\lambda))^j}{\lambda(h)^{k+j}}A_\lambda^{k-j}\ .
	\]
\end{Lem}

\begin{Rem}
	In passing, note that $b_{k-2,k}=b_{k-1,k}=\frac{(2k-2)!}{2^{k-1}(k-1)!}$. We remark also that $\theta_n(z)=\sum_{j=0}^nb_{j,n+1}z^{n-j}$ are so-called \emph{Bessel polynomials} \cite{grosswald-besselpol}, \cite[A001498]{OEIS}.
\end{Rem}

\begin{Par}
	Let $\lambda\in\bar\Sigma_1^+$, $\lambda(A_\lambda)=0$. By \thmref{Lem}{radialfmla}, we find \fa $I$, $\Abs0I=k$, that $\gamma_h(z_I\tilde z_I)=(-1)^{\frac12k(k+1)}\lambda(h)^{-k}A_\lambda^k$ ($h\in\ger a'$). Hence, 
	\[
		\bigcap\nolimits_{D\in\mathcal R_\lambda}\dom D=\bigcap\nolimits_{k=1}^{\frac12m_{1,\lambda}}\dom\lambda^{-k}\partial(A_\lambda)^k\ .
	\]
	
	The situation in the case $\lambda(A_\lambda)\neq0$ is different and requires a more detailed study.
\end{Par}
	
\begin{Par}[radial-action]
		Let $\lambda\in\bar\Sigma_1^+$, $\lambda(A_\lambda)\neq0$. Then $\cplxs[\ger a]\cong R[\lambda]$ where $R=\cplxs[\ker\lambda]$. This isomorphism is equivariant for $S(\cplxs A_\lambda)$ if we define an action $\partial$ on $R[\lambda]$ by requiring that $\partial(A_\lambda)$ be the unique $R$-derivation for which $\partial(A_\lambda)\lambda=\lambda(A_\lambda)$. 
		
		Now, let $R$ be an arbitrary commutative unital $\cplxs$-algebra. We define an action $\partial$ of $S(\cplxs A_\lambda)$ on $R[\lambda,\lambda^{-1}]$ by requiring that $\partial(A_\lambda)$ be the unique $R$-derivation \scth $\partial(A_\lambda)=\lambda(A_\lambda)$ and $\partial(A_\lambda)\lambda^{-1}=-\lambda(A_\lambda)\lambda^{-2}$. The action $\partial$ is faithful, because $\lambda(A_\lambda)\neq0$.
		
		Let $\mathcal D_\lambda$ be the subalgebra of $\End[_\cplxs]0{R[\lambda,\lambda^{-1}]}$ generated by $\partial(S(\cplxs A_\lambda))$ and $\cplxs[\lambda,\lambda^{-1}]$. In particular, we may embed $\mathcal R_\lambda\subset\mathcal D_\lambda$. We consider the action of $D\in\mathcal R_\lambda$, $D(h)=\gamma_h(z_I\tilde z_I)$, $\Abs0I=k$, on $p=\sum_{j=0}^Na_j\lambda^j\in R[\lambda]$, 
		\[
			Dp=(-1)^{\frac{k(k+1)}2}\textstyle\sum_{j=1}^Na_j\lambda(A_\lambda)^k\lambda^{j-2k}\sum_{i=(k-j)_+}^{k-1}(-1)^i(j)_{k-i}b_{ik}\in R[\lambda,\lambda^{-1}]\ .
		\]
		Since $\lambda(A_\lambda)\neq0$, we have $Dp\in R[\lambda]$ if and only if  
		\[
			a_j\textstyle\sum_{i=(k-j)_+}^{k-1}(-1)^i(j)_{k-i}b_{ik}=0\mathfa j=1,\dotsc,2k-1\ .
		\]
		
		We need to determine when the number
		\begin{equation}\label{eq:ajk-def}
		a_{jk}=\sum_{i=(k-j)_+}^{k-1}(-1)^i(j)_{k-i}b_{ik}=\sum_{i=(k-j)_+}^{k-1}\Parens2{-\frac12}^i(j)_{k-i}\frac{(k-1+i)!}{(k-1-i)!i!}
		\end{equation}
		is non-zero. 
\end{Par}
	
\begin{Par}
	Fix $k\sge1$. For $x\in\reals$ and $1\sle j\sle k$, let
	\[
		a_{2k-j,k}(x)=\sum_{i=0}^{k-1}x^i(2k-j)_{k-i}\frac{(k-1+i)!}{(k-1-i)!i!}\ .
	\]
	We claim that 
	\begin{equation}\label{eq:coeffx}
		a_{2k-j,k}(x)=\textstyle\frac{(j-1)!(2k-j)!}{(k-1)!}\sum_{\ell=0}^{j-1}\binom{k-1}\ell\binom{k-1}{j-1-\ell}x^\ell(1+x)^{k-1-\ell}\ .
	\end{equation}
	To that end, we rewrite 
	\[
		a_{2k-j,k}(x)=\frac{(j-1)!(2k-j)!}{(k-1)!}\sum_{i=0}^{k-1}\binom{k-1}i\binom{k+i-1}{j-1}x^i\ .
	\]
	Then, for fixed $x\in\reals$, we form the generating function
	\[
		f(z)=\sum_{j=1}^\infty z^{j-1}\sum_{i=0}^{k-1}\binom{k-1}i\binom{k+i-1}{j-1}x^i\ .
	\]
	It is easy to see 
	\begin{align*}
		f(z)&=\sum_{i=0}^{k-1}\binom{k-1}ix^i\sum_{j=1}^{k+i}\binom{k+i-1}{j-1}z^{j-1}\\
		&=(1+z)^{2k-2}\sum_{i=0}^{k-1}\binom{k-1}ix^i\Parens2{\frac1{1+z}}^{k-1-i}\\
		&=(1+z)^{k-1}((1+z)x+1)^{k-1}\ .
	\end{align*}
	
	On the other hand, we may form the generating function for the right hand side of (\ref{eq:coeffx}), 
	\[
		g(z)=\sum_{j=1}^\infty z^{j-1}\sum_{\ell=0}^{j-1}\binom{k-1}\ell\binom{k-1}{j-1-\ell}x^\ell(1+x)^{k-1-\ell}\ .
	\]
	Then
	\begin{align*}
		g(z)&=\sum_{\ell=0}^{k-1}\binom{k-1}\ell x^\ell(1+x)^{k-1-\ell}\sum_{j=\ell+1}^{k+\ell}\binom{k-1}{j-1-\ell}z^{j-1}\\
		&=\sum_{\ell=0}^{k-1}\binom{k-1}\ell(xz)^\ell(1+x)^{k-1-\ell}\sum_{j=0}^{k-1}\binom{k-1}jz^j\\
		&=(xz+x+1)^{k-1}(1+z)^{k-1}=f(z)\ .
	\end{align*}	
	Since the generating functions coincide, we have proved (\ref{eq:coeffx}). 
\end{Par}

\begin{Par}[evencoeffcase1]
	We notice that for $k\sge 1$ and $j=1,\dotsc,k$, $k-(2k-j)=j-k\sle0$, so $a_{2k-j,k}=a_{2k-j,k}\Parens1{-\tfrac12}$ by (\ref{eq:ajk-def}). By (\ref{eq:coeffx}), we obtain
	\[
		a_{2k-j,k}=\frac{(j-1)!(2k-j)!}{2^{k-1}(k-1)!}\sum_{\ell=0}^{j-1}(-1)^\ell\binom{k-1}\ell\binom{k-1}{j-1-\ell}
	\]
	For $j=1$, one gets
	\[
		a_{2k-1,k}=\frac{(2k-1)!}{2^{k-1}(k-1)!}\neq0\ .
	\]
	Now, let $j=2n$ where $1\sle n\sle\lfloor\tfrac k2\rfloor$. Then $\ell\mapsto(-1)^\ell\binom{k-1}\ell\binom{k-1}{2n-1-\ell}$ is odd under the permutation $\ell\mapsto2n-1-\ell$ of $\{0,\dotsc,2n-1\}$, so 
	\[
		a_{jk}=0\mathfa j=k,\dotsc,2k-2\ ,\ j\equiv0\ (2)\ .
	\]
\end{Par}

\begin{Par}[evencoeffcase2]
	Next, we study the behaviour of $a_{k-j,k}$ for $k\sge1$ and $j=1,\dotsc,k-1$, by a similar scheme. To that end, write 
	\begin{align*}
		a_{k-j,k}&=\sum_{i=j}^{k-1}\frac{(k-j)!(k-1+i)!}{(i-j)!(k-1-i)!i!}\Parens2{-\frac12}^i\\
		&=\frac{(k-1+j)!(k-j)!}{(k-1)!}\sum_{i=j}^{k-1}\binom{k-1}i\binom{k-1+i}{k-1+j}\Parens2{-\frac12}^i\ .
	\end{align*}
 	Observe that we may sum over $i=0,\dotsc,k-1$ since the second binomal coefficient vanishes for $i<j$. 
	
	Now, we fix $x\in\reals$ and define $f(z)=\sum_{j=1}^{k-1}a_{k-j,k}(x)z^{k+j-1}\in\cplxs[z]$ where
	\[
		a_{k-j,k}(x)=\sum_{i=0}^{k-1}\binom{k-1}i\binom{k-1+i}{k-1+j}x^i\ .
	\]
	We wish to study the coefficients of the polynomial $f$. Observe that the lowest power of $z$ occuring in $f(z)$ is $z^k$. Thus, we compute, modulo $\cplxs[z]_{<k}$, 
	\begin{align*}
		f(z)&=\sum_{i=0}^{k-1}\binom{k-1}ix^i\sum_{j=1}^i\binom{k-1+i}{k-1+j}z^{k+j-1}\\
		&=\sum_{i=0}^{k-1}\binom{k-1}ix^i\sum_{j=k}^{k-1+i}\binom{k-1+i}jz^j\\
		&\equiv(1+z)^{k-1}\sum_{i=0}^{k-1}\binom{k-1}i(x(1+z))^i=(1+z)^{k-1}(1+x(1+z))^{k-1}\ .
	\end{align*}	
	
	For $j=k,\dotsc,2k-2$, $a_{2k-j-1,k}(x)$ is the coefficient of $z^j$ in $f(z)$. Since 
	\[
		(1+z)^{k-1}(1+x(1+z))^{k-1}=\sum_{j=0}^{2k-2}z^j\sum_{i=0}^j\binom{k-1}{j-i}\binom{k-1}i(1+x)^{k-1-i}x^i\ ,
	\]
	we find, for $j=k,\dotsc,2k-2$,
	\begin{align*}
		a_{2k-j-1,k}(x)&=\sum_{i=0}^j\binom{k-1}{j-i}\binom{k-1}i(1+x)^{k-1-i}x^i\\
		&=(1+x)^{k-1}\sum_{i=j-k+1}^{k-1}\binom{k-1}{j-i}\binom{k-1}i\Parens2{\frac x{1+x}}^i\ .
	\end{align*}
	
	In particular,  
	\[
		a_{2k-j-1,k}\Parens1{-\tfrac12}=2^{1-k}\sum_{i=j-k+1}^{k-1}(-1)^i\binom{k-1}{j-i}\binom{k-1}i\ .
	\]
	Notice that the function $i\mapsto(-1)^i\binom{k-1}{j-i}\binom{k-1}i$ has parity $j$ with respect to the permutation $i\mapsto j-i$ of $\{j-k+1,\dotsc,k-1\}$. Since $2k-j-1$ is even and only if $j$ is odd, this implies
	\[
		a_{jk}=0\mathfa j=2,\dotsc,k-1\ ,\ j\equiv0\ (2)\ .
	\]
	
	We summarise the above considerations in the following proposition.
\end{Par}

\begin{Prop}[commondom]
	Let $R$ be a commutative unital $\cplxs$-algebra, and $\lambda\in\bar\Sigma_1^+$ \scth $\lambda(A_\lambda)\neq0$. Let $m\sge1$ be an integer, and for $k=1,\dotsc,m$, define
	\[
		D_k=(-1)^{\frac{k(k+1)}2}\textstyle\sum_{j=0}^{k-1}\frac{(k-1+j)!}{2^j(k-1-j)!j!}\frac{(-\lambda(A_\lambda))^j}{\lambda^{k+j}}A_\lambda^{k-j}\in \mathcal D_\lambda\ .
	\]
	Let $p=\sum_{j=0}^Na_j\lambda^j\in R[\lambda]$. Then $D_kp\in R[\lambda]$ \fa $k=1,\dotsc,m$ if and only 
	$a_j=0$ for all $j=1,\dotsc,2m-1$, $j\equiv1\ (2)$. 
\end{Prop}

\begin{proof}
	Let $1\sle k\sle m$. We have $a_{2k-1}a_{2k-1,k}=0$ and $a_{2k-1,k}\neq0$, so $a_{2k-1}=0$. Conversely, there are no further conditions, since $a_{km}=0$ for even $k$, $1<k<2m$. 
\end{proof}

\begin{Par}
	To apply \thmref{Prop}{commondom} to the determination of the image of the restriction map, let $\lambda\in\bar\Sigma_1^+$, $\lambda(A_\lambda)\neq0$. Note that $\cplxs[\ger a]=\cplxs[\ker\lambda][\lambda]$. Then \fa $p\in\cplxs[\ger a]$, 
	\[
		p=\textstyle\sum_{j=0}^\infty(j!)^{-1}\partial(A_\lambda)^jp|_{\ker\lambda}\Parens1{\frac\lambda{\lambda(A_\lambda)}}^j\ .
	\]
	\emph{I.e.}, if we take $R=\cplxs[\ker\lambda]$, then $p=\sum_ja_j\lambda^j$ where the coefficients are given by $a_j=\frac1{\lambda(A_\lambda)^jj!}\partial(A_\lambda)^jp|_{\ker\lambda}\in R$. Also, $\partial(A_\lambda)^ip|_{\ker\lambda}=0$ \fa $i=1,\dotsc,j$ if and only if $p\in\cplxs\oplus\lambda^{j+1}\cplxs[\ger a]$. Together with \thmref{Th}{super-chevalley}, we immediately obtain our main result, as follows. 
\end{Par}

\begin{Th}[chevrestrb]
	The restriction homomorphism $I(\ger p^*)\to S(\ger a^*)$ is a bijection onto the subspace $I(\ger a^*)=\bigcap_{\lambda\in\bar\Sigma_1^+}S(\ger a^*)^W\cap I_\lambda$ where
	\[
		I_\lambda=\textstyle\bigcap_{j=1}^{\frac12m_{1,\lambda}}\dom\lambda^{-j}\partial(A_\lambda)^j\mathtxt{if}\lambda(A_\lambda)=0
	\]
	and if $\lambda(A_\lambda)\neq0$, then $I_\lambda$ consists of those $p\in\cplxs[\ger a]$ \scth
	\[
		\partial(A_\lambda)^kp|_{\ker\lambda}=0\mathfa[odd integers] k\ ,\ 1\sle k\sle m_{1,\lambda}-1\ .
	\]
\end{Th}

\section{Examples}

\subsection{Scope of the theory}

\begin{Par}
	As remarked in \thmref{Par}{grouptype}, \thmref{Th}{chevrestrb} applies to a symmetric superpair of group type where $\ger k$ is classical and carries a non-degenerate invariant even form. The assumptions are still fulfilled if we add to $\ger k$ an even reductive ideal. Hence, $\ger k$ may be a direct sum of a reductive Lie algebra, and copies of any of the following Lie superalgebras \cite{kac-liesuperalgs}: 
	\begin{gather*}
		\ger{gl}(p|q,\cplxs)\ ,\ \sll(p|q,\cplxs)\ (p\neq q)\ ,\ \sll(p|p,\cplxs)/\cplxs\ ,\\
		\ger{osp}(p|2q,\cplxs)\ ,\ D(1,2;\alpha)\ ,\ F(4)\ ,\ G(3)\ .
	\end{gather*}
	As follows from \thmref{Prop}{superroot} (iv), in this situation one has $\lambda(A_\lambda)=0$ \fa $\lambda\in\bar\Sigma_1^+$.
\end{Par}

\begin{Par}[glexample]
	If we take $(\ger g,\ger k)$ to be an arbitrary reductive symmetric superpair, then the assumption of \emph{even type} amounts to an additional condition. 
	
	As an example, we consider $\ger g=\ger{gl}(p+q|r+s,\cplxs)$, $p,q,r,s\sge0$, where $\theta$ is given by conjugation with the diagonal matrix whose diagonal entries are the matrix blocks $1_p$, $-1_q$, $1_r$, $-1_s$. Let $\ger a\subset\ger p_0$ be the maximal Abelian subalgebra of all matrices
	\[
		\begin{Matrix}10&A&0&0\\-A^t&0&0&0\\0&0&0&B\\0&0&-B^t&0\end{Matrix}\in\cplxs^{(p+q+r+s)\times(p+q+r+s)}
	\]
	where $A=(D,0)$ or $A=\begin{Matrix}0D\\0\end{Matrix}$ for a diagonal matrix $D\in\cplxs^{\min(p,q)\times\min(p,q)}$, and similarly for $B$. Let $x_j$, $j=1,\dotsc,\min(p,q)$, and $y_\ell$, $\ell=1,\dotsc,\min(r,s)$, be the linear forms on $\ger a$ given by the entries of the diagonal blocks of $A,B$. 
	
	Consider the $\ger a$-module $\ger g_1$. Then the non-zero weights are 
	\[
		\pm(x_j\pm y_\ell)\ (2)\ ,\ \pm x_j\ (2\Abs0{r-s})\ ,\ \pm y_\ell\ (2\Abs0{p-q})
	\]
	with multiplicities given in parentheses \cite{schmittner-zirnbauer}. The sum $U\subset\ger g_1$ of the non-zero weight spaces therefore has dimension
	\begin{align*}
		8\min(p,q)\min(r,s)&+4\Abs0{r-s}\min(p,q)+4\Abs0{p-q}\min(r,s)\\
		&=2\Parens1{(p+q)(r+s)-\Abs0{p-q}\Abs0{r-s}}\ .
	\end{align*}
	(The equation follows by applying the formula $2\min(a,b)=a+b-\Abs0{a-b}$.)
	
	We have that $U$ is $\theta$-stable, and the action of a generic $h\in\ger a$ induces an automorphism of $U$. Hence, we have $\dim U_{\ger k}=\dim U_{\ger p}=\frac12\dim U$ where $U_{\ger k}$ and $U_{\ger p}$ are the projections of $U$ onto $\ger k_1$ and $\ger p_1$, respectively. It follows that $\dim U_{\ger p}=(p+q)(r+s)-\Abs0{p-q}\Abs0{r-s}$. On the other hand, 
	\[
		\dim\ger p_1=2(ps+rq)=(p+q)(r+s)-(p-q)(r-s)\ .
	\]
	Hence, $\ger z_{\ger p_1}(\ger a)=0$ if and only if $(p-q)(r-s)\sge0$, and $(\ger g,\ger k)$ is of even type  if and only if this condition holds. 
	
	We remark that in this case, the set $\bar\Sigma_1^+$ consists of the weights $x_j\pm y_\ell$ (for a suitably chosen positive system). For each $\lambda\in\bar\Sigma_1^+$, one has $\lambda(A_\lambda)=0$. 
\end{Par}

\begin{Par}
	A similar example arises by restricting the involution from \thmref{Par}{glexample} to the subalgebra $\ger g=\ger{osp}(p+q|r+s,\cplxs)$, where we now assume $r$ and $s$ to be even. We realise $\ger g$ by taking the direct sum of the standard non-degenerate symmetric forms on $\cplxs^p\oplus\cplxs^q$, and the direct sum of the standard symplectic forms on $\cplxs^r\oplus\cplxs^s$. 
	
	For $k$ even, denote by $J_k\in\cplxs^{k\times k}$ the matrix representing the standard symplectic form. Let $\ger a\subset\ger p_0$ be the maximal Abelian subalgebra of all matrices
	\[
		\begin{Matrix}10&A&0&0\\-A^t&0&0&0\\0&0&0&B\\0&0&J_sB^tJ_r&0\end{Matrix}\in\cplxs^{(p+q+r+s)\times(p+q+r+s)}
	\]
	where $A=(D,0)$ or $A=\begin{Matrix}0D\\0\end{Matrix}$ for a diagonal matrix $D\in\cplxs^{\min(p,q)\times\min(p,q)}$, and $B=(D',0)$ or $B=\begin{Matrix}0D'\\0\end{Matrix}$ for a diagonal matrix $D'\in\cplxs^{\frac12\min(r,s)\times\frac12\min(r,s)}$. 
	
	By restriction, we obtain the following non-zero $\ger a$-weights in $\ger g_1$, 
	\[
		\pm(x_j\pm y_\ell)\ (2)\ ,\ \pm x_j\ (\Abs0{r-s})\ ,\ \pm y_\ell\ (2\Abs0{p-q})\ ,
	\]
	where now $j=1,\dotsc,\min(p,q)\,,\,\ell=1,\dotsc,\tfrac12\min(r,s)$, and the multiplicities are given in parentheses \cite{schmittner-zirnbauer}.
	
	Let $U$ be the sum of all weight spaces for non-zero weights of the $\ger a$-module $\ger g_1$. Then the dimension of $U$ is 
	\begin{align*}
		4\min(p,q)\min(r,s)&+2\Abs0{r-s}\min(p,q)+2\Abs0{p-q}\min(r,s)\\
		&=(p+q)(r+s)-\Abs0{p-q}\Abs0{r-s}\ .
	\end{align*}
	If $U_{\ger p}$ is the projection of $U$ onto $\ger p_1$, then by the same argument as in \thmref{Par}{glexample}, $\dim U_{\ger p}=\frac12\dim U$. We have 
	\[
	\dim\ger p_1=pq+rs=\tfrac12\Parens1{(p+q)(r+s)-(p-q)(r-s)}\ ,
	\]
	so, as above, $(\ger g,\ger k)$ is of even type if and only if $(p-q)(r-s)\sge0$. In this case, as in \thmref{Par}{glexample}, the set $\bar\Sigma_1^+$ consists of the weights $x_j\pm y_\ell$ (for a suitable choice of positive system), and again we have $\lambda(A_\lambda)=0$ \fa $\lambda\in\bar\Sigma_1^+$. 
\end{Par}

\subsection{An extremal class: $\ger g=C(q+1)=\ger{osp}(2|2q,\cplxs)$, $\ger k_0=\ger{sp}(2q,\cplxs)$}

\begin{Par}
	Consider the Lie superalgebra $\ger g=C(q+1)=\ger{osp}(2|2q,\cplxs)$ where $q\sge1$ is arbitrary. Let $I=\begin{Matrix}00&1\\1&0\end{Matrix}\in\cplxs^{2\times2}$ and $J=\begin{Matrix}00&1\\-1&0\end{Matrix}\in\cplxs^{2q\times2q}$. If we realise $\ger g$ with respect to the orthosymplectic form $I\oplus J$, it consists of the matrices
	\[
		x=\begin{Matrix}1a&0&-w^{\prime t}&z^{\prime t}\\0&-a&-w^t&z^t\\ z&z'&A&B\\ w&w'&C&-A^t\end{Matrix}	
	\]
	where $a\in\cplxs$, $z,z',w,w'\in\cplxs^q$, $A,B=B^t,C=C^t\in\cplxs^{q\times q}$.

	The matrix $g=\begin{Matrix}0I&0\\0&1\end{Matrix}\in\cplxs^{(2+2q)\times(2+2q)}$ represents an even automorphism of the super-vector space $\cplxs^{2|2q}$, of order $2$. Since $g$ leaves the orthosymplectic form invariant, $\theta(x)=gxg$ defines an involutive automorphism of $\ger g$. Moreover, since $g^2=1$, the supertrace form $b(x,y)=\str(xy)$ on $\ger g$ is $\theta$-invariant. Hence, $(\ger g,\ger k)$, where $\ger k=\ger g_\theta$, is a reductive symmetric superpair.
	
	We compute
	\[
		\theta(x)=\begin{Matrix}1-a&0&-w^t&z^t\\0&a&-w^{\prime t}&z^{\prime t}\\z'&z&A&B\\w'&w&C&-A^t\end{Matrix}
	\]
	when $x\in\ger g$ is written as above. Hence, the general elements of $\ger k$ and $\ger p$ are respectively of the form 
	\[
		x=\begin{Matrix}10&0&-w^t&z^t\\0&0&-w^t&z^t\\z&z&A&B\\w&w&C&-A^t\end{Matrix}\nd
		x=\begin{Matrix}1a&0&w^t&-z^t\\0&-a&-w^t&z^t\\z&-z&0&0\\w&-w&0&0\end{Matrix}\ .
	\]
	
	It is immediate that the one-dimensional space $\ger a=\ger p_0$ is self-centralising in $\ger p_0$. In particular, any non-zero element of $\ger a$ is $b$-anisotropic (since $\ger p_0$ is non-degenerate). The bracket relation for the general element of $[\ger a,\ger g_1]$
	\[
		\left[\begin{Matrix}10&a&0&0\\-a&0&0&0\\0&0&0&0\\0&0&0&0\end{Matrix},\begin{Matrix}10&0&-w^{\prime t}&z^{\prime t}\\0&0&-w^t&z^t\\z&z'&0&0\\w&w'&0&0\end{Matrix}\right]=\begin{Matrix}10&0&-aw^{\prime t}&az^{\prime t}\\0&0&aw^t&-az^t\\-az&az'&0&0\\-aw&aw'&0&0\end{Matrix}
	\]
	implies in particular that $\ger z_{\ger p_1}(\ger a)=0$. Hence, $\ger a$ is an even Cartan subspace, and $(\ger g,\ger k)$ is of even type.
	
	Also, there are only two restricted roots, $\pm\lambda$, where $\lambda$ maps $x\in\ger a$ (as above) to $a$. Necessarily, $\lambda$ is odd, so $2\lambda\not\in\Sigma=\{\pm\lambda\}$, and $W=W(\Sigma_0)=1$. Since $A_\lambda$ is $b$-anisotropic, we have $\lambda(A_\lambda)\neq0$.
	
	Moreover, we must have $\ger p_1=\ger p_1^\lambda$, and this space has dimension $2q$, so $m_{1,\lambda}=2q$. From \thmref{Th}{chevrestrb}, we obtain the following result.
\end{Par}

\begin{Prop}[specialosp22q]
	Let $\ger g=\ger{osp}(2|2q,\cplxs)$, with the involution defined above. The image of the restriction map $S(\ger p^*)^{\ger k}\to S(\ger a^*)=\cplxs[\lambda]$ is 
	\[
		I(\ger a^*)=\Set1{p=\textstyle\sum_ja_j\lambda^j}{a_{2j-1}=0\ \forall j=1,\dotsc,q}\ .
	\]
	
	In particular, the algebra $I(\ger a^*)$ is isomorphic to the commutative unital $\cplxs$-algebra defined by the generators $\lambda_2$, $\lambda_{2q+1}$, and the relation
	\[
		(\lambda_2)^{2q+1}=(\lambda_{2q+1})^2\ .
	\]
\end{Prop}

\begin{proof}
	We only need to prove the presentation of $I(\ger a^*)$. Let $A$ be the unital commutative $\cplxs$-algebra defined by the above generators and relations. It is clear that there is a surjective algebra homomorphism from $\phi:A\to I(\ger a^*)$, defined by $\phi(\lambda_n)=\lambda^n$. 
		
	Consider on $I(\ger a^*)$ the grading induced by $\cplxs[\lambda]$. For any multiindex $\alpha=(\alpha_2,\alpha_{2q+1})$, define $\lambda_\alpha=(\lambda_2)^{\alpha_2}(\lambda_{2q+1})^{\alpha_{2q+1}}$ in the free commutative algebra $\cplxs[\lambda_2,\lambda_{2q+1}]$. The latter is graded via $\Abs0{\lambda_\alpha}=\Abs0\alpha=2\alpha_2+(2q+1)\alpha_{2q+1}$. The relation defining $A$ is homogeneous for this grading, so that $A$ inherits a grading. 

	By definition, $\phi$ respects the grading, and in fact, it is surjective in each degree of the induced filtration (and hence, in each degree of the grading). The relation of $A$ ensures that the image of $\lambda_\alpha$ in $A$, for any $\alpha$, depends only on $\Abs0\alpha$. Hence, $\dim A_j\sle1$ \fa $j$. This proves that $\phi$ is injective. 
\end{proof}

\begin{Cor}[singular]
	Under the assumptions of \thmref{Prop}{specialosp22q}, the algebra $I(\ger a^*)$ defines the singular curve in $\cplxs^2$ given by the equation $z^{2q+1}=w^2$. 
\end{Cor}

\begin{Par}
	We substantiate the above by some explicit computations. We have 
	\[
		\str\begin{Matrix}1a&0&w^t&-z^t\\0&-a&-w^t&z^t\\z&-z&0&0\\w&-w&0&0\end{Matrix}\begin{Matrix}1a'&0&w'^t&-z'^t\\0&-a'&-w'^t&z'^t\\z'&-z'&0&0\\w'&-w'&0&0\end{Matrix}=2aa'+4(w^tz'-z^tw')
	\]
	for the trace form $b$ on $\ger p=\ger a\oplus\ger p_1^\lambda$. In particular, 
	\[
		A_\lambda=\tfrac12\left(\begin{Matrix}[0]01&0&0&0\\0&-1&0&0\\0&0&0&0\\0&0&0&0\end{Matrix}\right)\ ,\ \lambda(A_\lambda)=\tfrac12\ .
	\]
	Setting 
	\begin{gather*}
		z_i=\tfrac12\left(\begin{Matrix}[0]00&0&0&-e_i^t\\0&0&0&e_i^t\\e_i&-e_i&0&0\\0&0&0&0\end{Matrix}\right)\ ,\ \tilde z_i=\tfrac12\left(\begin{Matrix}[0]00&0&e_i^t&0\\0&0&-e_i^t&0\\0&0&0&0\\e_i&-e_i&0&0\end{Matrix}\right)\ ,\\
		y_i=\tfrac12\left(\begin{Matrix}[0]00&0&0&-e_i^t\\0&0&0&-e_i^t\\-e_i&-e_i&0&0\\0&0&0&0\end{Matrix}\right)\ ,\ \tilde y_i=\tfrac12\left(\begin{Matrix}[0]00&0&e_i^t&0\\0&0&e_i^t&0\\0&0&0&0\\-e_i&-e_i&0&0\end{Matrix}\right)\ ,
	\end{gather*}
	one verifies the conditions from \thmref{Par}{symplbasis}, namely
	\begin{gather*}
		y_i,\tilde y_i\in\ger k_1\,,\,z_i,\tilde z_i\in\ger p_1\,,\,y_i+z_i,\tilde y_i+\tilde z_i\in\ger g_1^\lambda\,,\,b(y_i,\tilde y_j)=b(\tilde z_j,z_i)=\delta_{ij}\,,\\
		b(y_i,y_j)=b(\tilde y_i,\tilde y_j)=b(z_i,z_j)=b(\tilde z_i,\tilde z_j)=0\ .
	\end{gather*}
	Then one computes
	\begin{gather*}
		[y_i,z_j]=[\tilde y_i,\tilde z_j]=0\ ,\ [y_i,\tilde z_j]=-\delta_{ij}A_\lambda\ ,\ [\tilde y_i,z_j]=\delta_{ij}A_\lambda\ ,\\
		[A_\lambda,y_i]=\tfrac12z_i\ ,\ [A_\lambda,z_i]=\tfrac12y_i\ ,\ [A_\lambda,\tilde y_i]=\tfrac12\tilde z_i\ ,\ [A_\lambda,\tilde z_i]=\tfrac12\tilde y_i\ .
	\end{gather*}
	
	Let $\zeta_i$, $\tilde\zeta_i$, $i=1,\dotsc,q$, be the basis of $\ger p_1^*$, dual to $z_i$, $\tilde z_i$, $i=1,\dotsc,q$, so 
	\[
		\Cdual0{\tilde z_j}{\zeta_i}=-\Cdual0{z_i}{\tilde\zeta_j}=\delta_{ij}\ ,\ \Cdual0{z_j}{\zeta_i}=\Cdual0{\tilde z_i}{\tilde\zeta_j}=0\ .
	\]
	Then $\Cdual0z{\zeta_i}=b(z,z_i)$, $\Cdual0z{\tilde\zeta_i}=b(z,\tilde z_i)$, and one has
	\begin{gather*}
		\ad^*(y_i)\zeta_j=\ad^*(\tilde y_i)\tilde\zeta_j=0\ ,\ -\!\ad^*(y_i)\tilde\zeta_j=\ad^*(\tilde y_i)\zeta_j=\delta_{ij}\lambda\ ,\\
		\ad^*(y_i)\lambda=-\tfrac12\zeta_i\ ,\ \ad^*(\tilde y_i)\lambda=-\tfrac12\tilde\zeta_i\ .
	\end{gather*}
	Also, we observe $\Cdual0{z_I\tilde z_Jh^\nu}{\zeta_K\tilde\zeta_L\lambda^\mu}=\delta_{IL}\delta_{JK}\delta_{\nu\mu}(-1)^{\Abs0I\Abs0J}\nu!\lambda(h)^\nu$. 
	
	The preimages $p_2,p_{2q+1}$ of the generators $\lambda^2,\lambda^{2q+1}$ in $S(\ger p^*)^{\ger k}$ under the restriction map can be deduced from \thmref{Par}{radial-action}, because $\ger p=\ger a\oplus\ger p_1^\lambda$. Indeed, let $P=\phi(p_N)$ where $N=2$ or $N=2q+1$. By the formulae from \thmref{Par}{radial-action}, for $q\sge\Abs0I=k>0$ and $h\in\ger a'$, 
	\begin{align*}
		\textstyle\sum_{\nu=0}^\infty\tfrac1{\nu!}\Cdual0{z_I\tilde z_Jh^\nu}{p_N}&=P(z_I\tilde z_J;h)=(\partial_{\gamma_h(z_I\tilde z_J)}\lambda^N)(h)\\
		&=\delta_{IJ}(-1)^{\tfrac12k(k+1)}2^{-k}a_{Nk}\lambda(h)^{N-2k}
	\end{align*}
	where
	\[
		a_{Nk}=\textstyle\sum_{i=(k-N)_+}^{k-1}\Parens1{-\tfrac12}^i(N)_{k-i}\frac{(k-1+i)!}{(k-1-i)!i!}\ .
	\]
	Thus, 
	\[
		p_N=\lambda^N+\textstyle\sum_{k=1}^{\min(N,q)}(-1)^{\tfrac12k(k+3)}2^{-k}a_{Nk}\lambda^{N-2k}\sum_{\Abs0I=k}\zeta_I\tilde\zeta_I\ .
	\]
	
	When $N=2$ and $k\sge2$, then $a_{Nk}=0$ by \thmref{Par}{evencoeffcase1} and \thmref{Par}{evencoeffcase2}. On the other hand, $a_{21}=2$. Hence,
	\[
		p_2=\lambda^2+\textstyle\sum_{i=1}^q\zeta_i\tilde\zeta_i
	\]
	is the super-Laplacian, and 
	\[
		p_{2q+1}=\lambda^{2q+1}+\textstyle\sum_{k=1}^q(-1)^{\tfrac12k(k+3)}2^{-k}a_{2q+1,k}\lambda^{2(q-k)+1}\sum_{\Abs0I=k}\zeta_I\tilde\zeta_I\ .
	\]
	These elements are clearly subject to the relation $p_2^{2q+1}=p_{2q+1}^2$. 
	
	One readily checks
	\[
		\ad^*(y_i)p_2=-\lambda\zeta_i+\zeta_i\lambda=0\nd\ad^*(\tilde y_i)p_2=-\lambda\tilde\zeta_i+\lambda\tilde\zeta_i=0\ .
	\]
	In case $q=1$, one has $p_3=\lambda^3+\frac32\lambda\zeta_1\tilde\zeta_1$, and
	\begin{align*}
		\ad^*(y_1)p_3&=-\tfrac32\lambda^2\zeta_1-\tfrac32\lambda\zeta_1\ad^*(y_1)\tilde\zeta_1=-\tfrac32\lambda^2\zeta_1+\tfrac32\lambda\zeta_1\lambda=0\ ,\\
		\ad^*(\tilde y_1)p_3&=-\tfrac32\lambda^2\tilde\zeta_1+\tfrac32\lambda\ad^*(\tilde y_1)(\zeta_1)\tilde\zeta_1=-\tfrac32\lambda^2\tilde\zeta_1+\tfrac32\lambda^2\tilde\zeta_1=0\ .
	\end{align*}
	
	To verify the $\ger k_0$-invariance, let 
	\[
		x=\left(\begin{Matrix}[0]00&0&0&0\\0&0&0&0\\0&0&A&B\\0&0&C&-A^t\end{Matrix}\right)\in\ger k_0=\ger{sp}(2q,\cplxs)\ .
	\]
	Then 
	\[
		\ad^*(x)\zeta_i=\textstyle\sum_{j=1}^q\Parens1{A_{ji}\zeta_j+C_{ji}\tilde\zeta_j}\nd\ad^*(x)\tilde\zeta_i=\sum_{j=1}^q\Parens1{B_{ji}\zeta_j+A_{ji}\tilde\zeta_j}\ .
	\]
	This implies
	\[
		\ad^*(x)(\zeta_i\tilde\zeta_i)=\textstyle\sum_{j\neq i}\Parens1{C_{ji}\tilde\zeta_j\tilde\zeta_i-B_{ji}\zeta_i\zeta_j}\ .
	\]
	Since $B=B^t$, $C=C^t$, we deduce $\sum_{i=1}^q\ad^*(x)(\zeta_i\tilde\zeta_i)=0$. Since $\ger a=\ger z(\ger g_0)$ and thus $\ad^*(\ger k_0)\lambda=0$, this implies that $p_2$ (for general $q$) and $p_3$ (for $q=1$) are $\ger k$-invariant. 
\end{Par}

\bibliographystyle{alpha}%
\bibliography{ahz-chevalley}%

\end{document}